\documentclass[11pt,reqno]{amsart}

\usepackage{amsmath,amssymb}

\usepackage[usenames,dvipsnames]{color}
\usepackage{hyperref}
\hypersetup{
    colorlinks=true,			
    linkcolor=blue,			
    citecolor=Magenta,		
    filecolor=magenta,		
    urlcolor=cyan			
}

\usepackage[all]{xy}

\newtheorem*{thm}{Theorem}
\newtheorem*{prop}{Proposition}
\newtheorem*{lem}{Lemma}
\newtheorem*{cor}{Corollary}

\newenvironment{pf}{\paragraph{{\sc Proof}}}{\qed\par\medskip}
\theoremstyle{definition}

\newtheorem*{rem}{Remark}
\newtheorem*{rems}{Remarks}

\numberwithin{equation}{section}
\numberwithin{figure}{section}

\newcommand{\cfun}[2]{S\lp\M_{#1,#2}^*\rp}
\newcommand{\qfun}[2]{S_{\hbar}\lp\M_{#1,#2}^*\rp}
\newcommand{\id}{\mathbf{1}}
\newcommand{\pde}[2]{\frac{\partial #1}{\partial #2}}
\newcommand{\KZ}{\operatorname{KZ}}
\newcommand{\el}[2]{E_{#1}^{(#2)}}
\newcommand{\ve}{\varepsilon}

\newcommand{\lp}{\left(}
\newcommand{\rp}{\right)}
\newcommand{\la}{\left\langle}
\newcommand{\ra}{\right\rangle}

\newcommand{\g}{\mathfrak{g}}
\newcommand{\h}{\mathfrak{h}}
\newcommand{\gl}{\mathfrak{gl}}
\newcommand{\Lsl}{\mathfrak{sl}}

\newcommand{\Sym}{\mathfrak{S}}

\newcommand{\fp}{\mathsf{h}}

\newcommand{\bfA}{\mathbf{A}}

\newcommand{\A}{\mathcal{A}}

\newcommand{\E}{\mathcal{E}}
\newcommand{\F}{\mathcal{F}}

\newcommand{\HH}{\mathcal{H}}

\newcommand{\J}{\mathcal{J}}
\newcommand{\K}{\mathcal{K}}
\newcommand{\LL}{\mathcal{L}}
\newcommand{\M}{\mathcal{M}}
\newcommand{\calN}{\mathcal{N}}

\newcommand{\RR}{\mathcal{R}}

\newcommand{\V}{\mathcal{V}}
\newcommand{\X}{\mathcal X}
\newcommand{\cZ}{\mathcal Z}

\newcommand{\C}{\mathbb{C}}
\newcommand{\nC}{\mathbb{C}^{\times}}
\newcommand{\N}{\mathbb{N}}

\newcommand{\Z}{\mathbb{Z}}
\newcommand{\qS}{\mathbb{S}}
\newcommand{\qL}{\mathbb{L}}

\newcommand{\fs}{{\mathbf s}}

\newcommand {\wh}[1]{\widehat{#1}}
\newcommand {\ol}[1]{\overline{#1}}
\newcommand {\ul}[1]{\underline{#1}}
\newcommand {\aff}{\scriptscriptstyle{\operatorname{aff}}}

\newcommand{\reg}{\scriptscriptstyle{\operatorname{reg}}}
\newcommand{\areg}{\scriptscriptstyle{\operatorname{a-reg}}}
\newcommand {\ddeg}{\operatorname{deg}}
\newcommand {\Aut}{\operatorname{Aut}}

\newcommand{\End}{\operatorname{End}}
\newcommand{\Hom}{\operatorname{Hom}}
\newcommand{\ad}{\operatorname{ad}}
\newcommand{\Ad}{\operatorname{Ad}}
\newcommand{\ev}{\operatorname{ev}}
\newcommand {\aand}{\qquad\text{and}\qquad}

\newcommand {\fd}{finite--dimensional }
\newcommand {\lhs}{left--hand side }
\newcommand {\rhs}{right--hand side }
\newcommand {\wrt}{with respect to }
\newcommand {\hw}{highest weight }

\newcommand{\ds}{\displaystyle}
\newcommand{\wt}[1]{\widetilde{#1}}

\newcommand{\zed}{\mathfrak{z}}

\newcommand{\qdet}{\operatorname{qdet}}
\newcommand{\Th}{\Theta}

\newcommand {\hloop}{U_\hbar(L\g)}
\newcommand {\hloopgl}[1]{U_\hbar(L\gl_{#1})}
\newcommand {\hloopsl}[1]{U_{\hbar}(L\Lsl_{#1})}
\newcommand {\Uhg}{U_\hbar\g}
\newcommand{\Uhgl}[1]{U_{\hbar}\gl_{#1}}
\newcommand{\Uhsl}[1]{U_{\hbar}\Lsl_{#1}}

\newcommand {\Yhgl}[1]{Y_\hbar(\gl_{#1})}

\newcommand {\Ysl}[1]{Y_{\fp}\Lsl_{#1}}

\newcommand {\isom}{\stackrel{\sim}{\rightarrow}}
\newcommand{\cbin}[2]{\left(\begin{array}{c} #1\\ #2\end{array}\right)}
\newcommand{\qbin}[3]{\left[\begin{array}{c} #1 \\ #2\end{array}\right]_{#3}}

\renewcommand {\tt}{\mathfrak t}
\renewcommand {\gl}{\mathfrak{gl}}
\renewcommand {\sl}{\mathfrak{sl}}

\newcommand {\Yh}{Y_\fp}
\newcommand {\IV}{\mathbb V}

\newcommand {\bfV}{V}
\newcommand {\bfcV}{\mathbf{\mathcal{V}}}
\newcommand {\KKZ}{_{KZ}}
\newcommand{\ctensor}{\wh{\otimes}}

\newcommand {\trigosl}[1]{\wh{\nabla}_C^{\sl_{#1}}}
\newcommand {\trigogl}[1]{\wh{\nabla}_C^{\gl_{#1}}}
\newcommand {\trigoslev}[2]{\wh{\nabla}_{C,\ul{#2}}^{\sl_{#1}}}
\newcommand {\trigoglev}[2]{\wh{\nabla}_{C,\ul{#2}}^{\gl_{#1}}}
\newcommand {\pisl}{\pi^{\sl_2}_{C,\ul{s}}}
\newcommand {\pigl}{\pi^{\gl_2}_{C,\ul{s}}}
\newcommand {\Ygl}[1]{Y_{\fp}\gl_{#1}}
\newcommand {\sfA}{{\mathsf A}}

\newcommand {\KM}{Kac--Moody }

\newcommand {\eg}{{\it e.g., }}

\newcommand {\qt}{quasitriangular }
\newcommand {\Id}{\operatorname{id}}
\renewcommand \qfun[2]{\C_\hbar[\M_{#1,#2}]}
\renewcommand \cfun[2]{\C[\M_{#1,#2}]}
\renewcommand {\ast}{\theta}

\renewcommand {\fs}{s}
\newcommand {\BSL}{B_{SL_2}}
\newcommand {\Bsl}{B_{\sl_2}}
\newcommand {\BGL}{B_{GL_2}}
\newcommand {\CP}{Chari--Pressley }
\newcommand {\EK}{Etingof--Kazhdan }
\newcommand {\EGS}{Etingof--Geer--Schiffmann }
\newcommand {\KD}{Kohno--Drinfeld }
\newcommand {\CYBE}{classical Yang--Baxter equations }
\newcommand {\fml}{[[\hbar]]}
\newcommand {\ind}{\operatorname{ind}}
\newcommand {\wtR}{R} 

\newcommand {\Uh}{U_\hbar}
\renewcommand {\ll}{\mathfrak l}

\newcommand {\half}[1]{\frac{#1}{2}}
\newcommand {\El}[1]{E_{#1}}
\newcommand {\XX}{\mathfrak{X}}

\xyoption{arc}
\begin{document}

\title{Monodromy of the trigonometric Casimir connection for $\Lsl_2$}
\author[S. Gautam]{Sachin Gautam}
\address{Mathematics Department,
Columbia University,
2990 Broadway,
New York, NY 10027}
\email{sachin@math.columbia.edu}
\address{Department of Mathematics,
Northeastern University,
360 Huntington Avenue,
Boston, MA 02115.}
\email{gautam.s@husky.neu.edu}
\author[V. Toledano Laredo]{Valerio Toledano Laredo}
\address{Department of Mathematics,
Northeastern University,
360 Huntington Avenue,
Boston, MA 02115.}
\email{V.ToledanoLaredo@neu.edu}
\thanks{Both authors are supported by NSF grants DMS--0707212 and DMS--0854792}
\begin{abstract}
We show that the monodromy of the trigonometric Casimir connection on the
tensor product of evaluation modules of the Yangian $\Ysl{2}$ is described by
the quantum Weyl group operators of the quantum loop algebra $\hloopsl{2}$.
The proof is patterned on the second author's computation of the monodromy
of the rational Casimir connection for $\Lsl_n$ via the dual pair $(\gl_k,\gl_n)$,
and rests ultimately on the \EGS computation of the monodromy of the trigonometric
KZ equations. It relies on two new ingredients: an affine extension of the duality
between the $R$--matrix of $\Uhsl{k}$ and the quantum Weyl group element
of $\Uhsl{2}$, and a formula expressing the quantum Weyl group action
of the coroot lattice of $SL_2$ in terms of the commuting generators of
$\hloopsl{2}$. Using this formula, we define quantum Weyl group operators
for the quantum loop algebra $\hloopgl{2}$, and show that they describe the
monodromy of the trigonometric Casimir connection on a tensor product of
evaluation modules of the Yangian $\Ygl{2}$.
\end{abstract}
\maketitle

\setcounter{tocdepth}{1}
\tableofcontents

\section{Introduction}

\subsection{} 

Let $\g$ be a complex, semisimple Lie algebra, $G$ the corresponding
connected and simply--connected Lie group, $H\subset G$ a maximal
torus and $W$ the corresponding Weyl group. In \cite{valerio3}, a flat
$W$--equivariant connection $\wh{\nabla}_C$ was constructed on $H$
which has logarithmic singularities on the root subtori of $H$ and values
in any \fd representation of the Yangian $Y_{\fp}\g$. By analogy with the
description of the monodromy of the rational Casimir connection obtained
in \cite{valerio4, valerio0}, it was conjectured in \cite{valerio3} that the
monodromy of the trigonometric Casimir connection $\wh{\nabla}_C$
is described by the action of the affine braid group $B_G$ of $G$ arising from
the quantum Weyl group operators of the quantum loop algebra $\hloop$.

\subsection{} 

The aim of the present paper is to prove this conjecture when $\g =\sl_2$
and $V$ is a tensor product of evaluation modules. Note that, by a theorem
of \CP \cite{chari-pressley-yangian}, such representations include all irreducible
$\Ysl{2}$--modules. To state our main result, let $V_1,\ldots,V_k$ be \fd
$\sl_2$--modules, $z_1,\ldots,z_k$ points in $\C$, and
\[\bfV(\ul{z})=V_1(z_1)\otimes\cdots\otimes V_k(z_k)\]
the tensor product of the corresponding evaluation representations of $
\Ysl{2}$. The monodromy of the trigonometric Casimir connection yields
an action of the affine braid group $\BSL$ on $\bfV(\ul{z})$.

Let $\V_i$ be a quantum deformation of $V_i$, that is a module over the
quantum group $U_\hbar\sl_2$ such that $\V_i/\hbar\V_i\cong V_i$. Set
$\hbar=4\pi\imath\fp$ and $\zeta_i=\exp(-\hbar z_a)$, and consider the
tensor product of evaluation representations of the quantum loop algebra
$U_\hbar(L\sl_2)$ given by
\[\bfcV (\ul{\zeta})= \V_1(\zeta_1)\otimes\cdots\otimes\V_k(\zeta_k)\]
The quantum Weyl group operators $\qS_0, \qS_1$ of $\hloopsl{2}$ yield a
representation of $\BSL$ on $\bfcV(\ul{\zeta})$ \cite{kirillov-reshetikhin,lusztig-book,
soibelman}. The main result of this paper is the following

\begin{thm}\label{thm: intro-theorem-1}
The monodromy action of the affine braid group $\BSL$ on $\bfV(\ul{z})$
is equivalent to  its quantum Weyl group action on $\bfcV(\ul{\zeta})$.
\end{thm}

\subsection{}\label{ss:intro-duality}

The proof of the above theorem relies on two dualities between the
Lie algebras $\sl_k$ and $\sl_n$ discovered in \cite{valerio4}\footnote
{in the case relevant to the present paper, $n=2$.}. The first duality
arises from their joint action on the space $\cfun{k}{n}$ of functions
on $k\times n$ matrices, and identifies the rational Casimir connection
of $\sl_k$ with the rational KZ connection on $n$ points for $\sl_k$. 
The second duality arises from the action of the corresponding
quantum groups $\Uhsl{k}$ and $\Uhsl{n}$ on a noncommutative deformation
of $\cfun{k}{n}$, and identifies the quantum Weyl group elements 
of $\Uhsl{n}$ with the $R$--matrices of $\Uhsl{k}$.

These dualities were used in \cite{valerio4} together with the \KD
theorem for $\sl_k$, to show that the monodromy of the rational
Casimir connection of $\sl_n$ is described by the quantum Weyl
group operators of $\Uhsl{n}$.

\subsection{} 

In this paper, we apply a similar strategy to compute the monodromy
of the trigonometric Casimir connection of $\sl_2$ and, in fact, $\gl_2$.
The latter connection is an extension of the former to the maximal torus
of $GL_2$ constructed in \cite{valerio3}, and takes values in the Yangian
$\Ygl{2}$. Its evaluation on a tensor product of evaluation modules 
coincides, up to abelian terms, with the trigonometric dynamical differential
equations considered in \cite{tarasov-varchenko-duality}. In particular,
we also compute the monodromy of these equations.

The duality between the Casimir and KZ connections identifies the
trigonometric Casimir connection of $\gl_2$ with the trigonometric
KZ connection of $\gl_k$ (see, \eg \cite{tarasov-varchenko-duality}).
In turn, the monodromy of the latter was computed by \EGS in terms
of data coming from the quantum group $\Uhgl{k}$ \cite{etingof-geer}.
This reduces the original problem to interpreting this data in terms of
the quantum loop algebra $\hloopgl{2}$.

Part of this interpretation, namely the one pertaining to the data
describing the monodromy of the finite braid group $\Z\cong\Bsl\subset
\BSL$, is provided by the duality between $\Uhgl{k}$ and $\Uhgl{2}$
of \cite{valerio4} alluded to in \ref{ss:intro-duality}. What remains is the
description of the operators giving the action of the
coroot lattice $\Z^2\cong Q^\vee\subset\BGL$ of $GL_2$, in
terms of appropriate, commuting quantum Weyl group operators
of $\hloopgl{2}$.

\subsection{}

To the best of our knowledge, quantum Weyl group operators
giving an action of the coroot lattice of $GL_2$ on \fd representations
of the quantum loop algebra $\hloopgl{2}$ have not been defined. 
Moreover, for $\hloopsl{2}$, no compact, explicit formula appears
to be known for the element $\qS_0\qS_1$ giving the action of the
generator of the coroot lattice of $SL_2$. In this paper, we give the
following solution to both of these problems.

Let $\tt\subset\gl_2$ and $\h\subset\sl_2$ be the Cartan subalgebras
of diagonal and traceless diagonal matrices respectively, and
$U_0\subset\hloopgl{2}$, $U_0'\subset\hloopsl{2}$ the commutative
subalgebras deforming $U(\tt[z,z^{-1}])$ and $U(\h[z,z^{-1}])$.
Then, we prove the following.

\begin{thm}\hfill
\begin{enumerate}
\item There exist elements $\qL_1, \qL_2$ in a completion of $U_0$
such that $\{\qS=\qS_1, \qL_1, \qL_2\}$ satisfy the defining relations
of the affine braid group $\BGL$.
\item The element $\qL=\qL_1\qL_2^{-1}$ lies in a completion of $U_0'$,
and coincides with the quantum Weyl group element $\qS_0\qS_1$
giving the action of the generator of the coroot lattice of $SL_2$.
\end{enumerate}
\end{thm}

The elements $\qL_1,\qL_2$ are given by explicit formulae in terms
of the generators of $U_0$. For $\qL=\qL_1\qL_2^{-1}$, these are as
follows. Let $\{H_k\}_{k\in\Z}$ be the generators of $U_0'$ with classical
limit $\{h\otimes z^k\}$, where $h$ is the standard generator of $\h$
(see Section \ref{sec: qla}). Define, for any $r\in\N$,
\[\wt{H}_r=
H_0\phantom{_{i,}}+ \sum_{s=1}^r (-1)^s \cbin{r}{s}\frac{s}{[s]} H_s\]
and note that $\wt{H}_r=h\otimes (1-z)^r$ mod $\hbar$.
Then, we show that
\[ \qL=\exp\lp \sum_{r\geq 1} \frac{\wt{H}_r}{r} \rp \]
thus extending to the $q$--setting the fact that the classical limit of
$\qL$ is the loop
\[z\mapsto
\begin{pmatrix}z^{-1}&0\\0&z\end{pmatrix}=
\exp(-h\log z)\]
The operators $\qL_1,\qL_2$ are given by similar formulae. These
generalise in fact to any complex semisimple Lie algebra and to 
$\gl_n$ \cite{sachin-valerio-4}.

\subsection{}

Once the operators $\qL_1,\qL_2$ are explicitly defined, a direct 
computation shows that their action on quantum $k\times 2$ matrix
space coincides with that of the $\Uhgl{k}$ operators which, by
\cite{etingof-geer} describe the monodromy of the trigonometric KZ
connection of $\gl_k$, thus providing an extension of the $q$--duality
of \cite{valerio4} to the affine setting. Theorem \ref{thm: intro-theorem-1},
and its analogue for $\gl_2$ follow as a direct consequence.

\subsection{} 

The results of the present paper extend without essential modification
to the case of $\g=\sl_n$ and $\gl_n$, and give a computation of the
monodromy of the trigonometric Casimir connection of $\g$ with
values in a tensor product of arbitrary \fd evaluation representations
of the Yangian $Y_\fp\g$, in terms of the quantum Weyl group
operators of the quantum loop algebra $\hloop$.

\subsection{Outline of the paper}

Sections \ref{sec: trig-casimir} and \ref{sec: gl-sl} review the
definition of the Yangian and trigonometric Casimir connections
of the Lie algebras $\sl_2$ and $\gl_2$ respectively. Section \ref
{se:braid gps} gives presentations of the affine braid groups
$\BSL$ and $\BGL$, and describes the embedding $\BSL\subset
\BGL$ resulting from the inclusion of the maximal tori of $SL_2$
and $GL_2$ in terms of the corresponding generators.

In Section \ref{se:trigo KZ}, we review the definition of the
trigonometric KZ connection for the Lie algebra $\gl_k$ and, in
Section \ref{sec: howe} the fact that, under $(\gl_k,\gl_2)$--duality,
the trigonometric Casimir connection for $\gl_2$ is identified with
the trigonometric KZ connection for $\gl_k$. In Section \ref{sec: etingof-geer}
we describe, following \cite{etingof-geer}, the monodromy of the
latter connection in terms of the quantum group $\Uhgl{k}$.

In Section \ref{sec: qla}, we review the definition of the quantum
loop algebras $\hloopgl{2}$ and $\hloopsl{2}$. Section \ref{sec: qWeyl-gl}
contains the main construction of this paper. We first extend the quantum
Weyl group action of the affine braid group $\BSL$ on $\hloopsl{2}$ to
one of $\BGL$ on $\hloopgl{2}$. We then show that this action is
essentially inner, by exhibiting elements in an appropriate completion
of the maximal commutative subalgebra of $\hloopgl{2}$, whose adjoint
action coincides with the quantum Weyl group action of the coroot
lattice of $\gl_2$.

Section \ref{sec: q-howe} describes the joint action of $\Uhgl{k}$
and $\Uhgl{2}$ on the space $\qfun{k}{2}$ of quantum $k\times 2$
matrices. In Section \ref{sec: equivalence}, we prove the equality of
two actions of the affine braid group $\BGL$ on $\qfun{k}{2}$. The
first arises from its structure as $\Uhgl{k}$--module, and describes
the monodromy of the trigonometric KZ equations; the second from
its structure as a tensor product of $k$ evaluation modules of
$\hloopgl{2}$.

In Section \ref{se:monodromy}, we prove that the monodromy
of the trigonometric Casimir connection for $\g=\sl_2$ (resp. $\g=\gl_2$)
on a tensor product of evaluation modules is described by the
quantum Weyl group operators of $\hloop$.

Appendix \ref{app-sec: etingof-geer} outlines the computation
of the monodromy of the trigonometric KZ connection given in
\cite{etingof-geer}. Appendix \ref{app:J^n} contains the proof
of a technical result bearing upon the completions of the quantum
loop algebras $\hloopsl{2}$ and $\hloopgl{2}$ required to handle
quantum Weyl group elements.

\subsection*{Acknowledgments} The present paper was completed 
while both authors visited the Kavli Institute for Theoretical Physics
at the University of California, Santa Barbara. We are grateful to the
organisers of the program {\it Nonperturbative Effects and Dualities
in QFT and Integrable Systems} for their invitation, and Edward
Frenkel and KITP for supporting our stay through their DARPA
and NSF grants HR0011-09-1-0015 and PHY05-51164 respectively.

\section{The trigonometric Casimir connection of $\sl_2$}\label{sec: trig-casimir}

\subsection{The Yangian $\Ysl{2}$ \cite{drinfeld-yangian-qaffine}}\label{sec: defn-yangian}

The Yangian $\Ysl{2}$ is the unital, associative algebra over $\C
[\fp]$ generated by elements $\{\xi_r, e_r, f_r\}_{r\in \N}$, subject
to the relations
\begin{itemize}
\item[(Y1)] For each $r,s\in \N$, \[[\xi_r, \xi_s]=0\]
\item[(Y2)] For each $r\in \N$,
\[[\xi_0, e_r]=2e_r \aand [\xi_0, f_r]=-2f_r\]
\item[(Y3)] For each $r,s\in \N$,
\[[e_r, f_s]=\xi_{r+s}\]
\item[(Y4)] For each $r,s\in \N$,
\begin{align*}
[\xi_{r+1}, e_s] - [\xi_r, e_{s+1}]&=\phantom{-}\fp\lp \xi_re_s + e_s\xi_r\rp\\
[\xi_{r+1}, f_s] - [\xi_r, f_{s+1}]	&=-\fp\lp \xi_rf_s + f_s\xi_r\rp
\end{align*}
\item[(Y5)] For each $r,s\in \N$,
\begin{align*}
[e_{r+1}, e_s] - [e_r, e_{s+1}]	&=\phantom{-}\fp\lp e_re_s + e_se_r\rp \\
[f_{r+1}, f_s] - [f_r, f_{s+1}]	&=-\fp \lp f_rf_s+f_sf_r\rp
\end{align*}
\end{itemize}

$\Ysl{2}$ is an $\N$--graded algebra with $\ddeg(x_r)=r$ and 
$\ddeg{\fp}=1$. Moreover, it is a Hopf algebra with coproduct
determined by
\[\Delta(x_0)=x_0\otimes 1 + 1\otimes x_0\]
for $x=e,f,\xi$, and
\begin{equation}\label{eq: delta-xi}
\Delta(\xi_1)=\xi_1\otimes 1+1\otimes \xi_1+\fp\lp\xi_0\otimes \xi_0
-2f_0\otimes e_0\rp
\end{equation}
Let $\{e,f,h\}$ be the standard basis of the Lie algebra $\sl_2$.
Then, the map
\[e\to e_0\qquad f\to f_0\qquad h\to\xi_0\]
defines an embedding of $\sl_2$ into $\Ysl{2}$. In particular,
$\Ysl{2}$ is acted upon by $\sl_2$ via the adjoint action. This
action is integrable since the graded components of $\Ysl{2}$
are finite--dimensional.

\subsection{The trigonometric Casimir connection \cite{valerio3}}
\label{ss:trigo sl}

Let $G=SL_2(\C)$, $H\subset G$ the maximal torus consisting
of diagonal matrices and $\h\subset\sl_2$ its Lie algebra. The
Weyl group $W\cong\Z_2$ of $G$ acts on $H$ and on the
centraliser of $\h$ in $\Yh\sl_2$.

The trigonometric Casimir connection of $\sl_2$ is the flat, $W
$--equivariant connection on $H$ with values in $\Yh\sl_2$ given
by
\[\trigosl{2}=d-\left(\frac{\fp\kappa}{e^{\alpha}-1}-t_1\right)d\alpha\]
where $\kappa=e_0f_0 + f_0e_0$ is the truncated Casimir element
of $\sl_2$, $\alpha\in\h^*$ is defined by $\alpha(h)=2$, $d\alpha$
is the corresponding translation--invariant one--form on $H$, and
\[t_1=\ds \xi_1 - \frac{\fp}{2}\xi_0^2\]

\subsection{Evaluation homomorphism}\label{sec: ev-yangian}

For any $s\in\C[\fp]$, there is an algebra homomorphism $\ev_
s:\Ysl{2} \to U\Lsl_2[\fp]$ which is equal to the identity on $\sl_2
\subset\Ysl {2}$ and is otherwise determined by \cite[Prop. 2.5]
{chari-pressley-yangian}
\begin{gather*}
t_1 \mapsto sh-\frac{\fp}{2}\kappa
\end{gather*}
Note that if $s\in\fp\C[\fp]$, $\ev_s$ maps elements of positive
degree in $\Ysl{2}$ to $\fp U\sl_2[\fp]$ and therefore extends
to a homomorphism $\wh{\Ysl{2}}\to U\sl_2[[\fp]]$, where $\wh
{\Ysl{2}}$ is the completion of $\Ysl{2}$ \wrt its grading.

\subsection{}\label{sec: trig-casimir-ev}

Let $k\in\N^*$, $\ul{s}=(s_1,\ldots, s_k)\in\C[\fp]^k$ and consider
the homomorphism
\[\ev_{\ul{s}}=\ev_{s_1}\otimes\cdots\otimes \ev_{s_k}\circ
\Delta^{(k)} : \Ysl{2}\to U\Lsl_2^{\otimes k}[\fp]\]
where $\Delta^{(k)} : \Ysl{2}\to \Ysl{2}^{\otimes k}$ is the 
iterated coproduct.

\begin{prop}\label{prop: trig-casimir-ev}
The image of $\trigosl{2}$ under the homomorphism $\ev_{\ul{s}}$
is the $U\Lsl_2^{\otimes k}[\fp]$--valued connection on $H$ given
by
\[\trigoslev{2}{s}=d - \left(\fp\frac{\Delta^{(k)}
(\kappa)}{e^{\alpha}-1}-A\right)d\alpha\]
where
\[A=
\sum_{a=1}^k s_a h^{(a)}-
\frac{\fp}{2} \sum_{a=1}^k\kappa^{(a)}-
2\fp\negthickspace\negthickspace\sum_{1\leq a<b\leq k}
\negthickspace\negthickspace f^{(a)}e^{(b)}\\[1.1 ex]\]
and, for any $x\in U\Lsl_2$, $x^{(a)}=1^{\otimes(a-1)}\otimes x\otimes
1^{\otimes(k-a)}$.
\end{prop}
\begin{pf}
The element $t_1$ defined in \S \ref{ss:trigo sl} satisfies $\Delta(t_1)=
t_1\otimes 1+1\otimes t_1-2\fp f_0\otimes e_0$, so that 
\[\Delta^{(k)}(t_1)=\sum_a t_1^{(a)}-2\fp\sum_{a<b}f_0^{(a)}e_0^{(b)}\]
The result follows since $\ev_s(t_1)=sh - \frac{\fp }{2}\kappa$.
\end{pf}

\section{The trigonometric Casimir connection of $\gl_2$}\label{sec: gl-sl}

\subsection{The Yangian $\Ygl{2}$ \cite{drinfeld-qybe}}\label{sec: yangian-gl}

$\Ygl{2}$ is the unital, associative algebra over $\C[\fp]$
generated by elements $\{t_{ij}^{(r)}\}_{1\leq i,j\leq 2, r
\geq 1}$, subject to the relations\footnote{we follow the
sign conventions of \cite{molev-yangian}.}
\[[t_{ij}^{(r+1)}, t_{kl}^{(s)}] - [t_{ij}^{(r)}, t_{kl}^{(s+1)}]=\fp 
\lp t_{kj}^{(r)} t_{il}^{(s)} - t_{kj}^{(s)}t_{il}^{(r)}\rp\]
for any $r,s\geq 0$, where  $t_{ij}^{(0)}=\fp^{-1}\delta_{ij}$.
These imply that $E_{ij}\mapsto t_{ij}^{(1)}$ gives an embedding
of $\gl_2$ into $\Ygl{2}$, and that $\Ygl{2}$ is an $\N$--graded
algebra with
\[\ddeg(t_{ij}^{(r)})=r-1 \aand \ddeg(\fp)=1\]
Moreover, $\Ygl{2}$ is a Hopf algebra with coproduct given by
\[\Delta(t_{ij}(u))=\sum_k t_{ik}(u)\otimes t_{kj}(u)\]
where $t_{ij}(u)=\fp\sum_{r\geq 0}t_{ij}^{(r)}u^{-r}$.

\subsection{The embedding $\Ysl{2}\subset\Ygl{2}$\cite{brundan-kleshev,molev-yangian}}\label{ss:embed}

Let $e(u),f(u),\xi(u)\in\Ysl{2}[[u^{-1}]]$ be the generating series
\[e(u)=\fp\sum_{r\geq 0}e_r u^{-r-1}\qquad
f(u)=\fp\sum_{r\geq 0}f_r u^{-r-1}\qquad
\xi(u)=1+\fp\sum_{r\geq 0}\xi_r u^{-r-1}\]
Then, the following defines an embedding of graded Hopf algebras
$\imath:\Ysl{2}\to\Ygl{2}$ \cite[Rem. 3.1.8]{molev-yangian}
\begin{gather*}
e(u)\mapsto t_{21}(u)t_{11}(u)^{-1}
\qquad\qquad\qquad
f(u)\mapsto t_{11}(u)^{-1}t_{12}(u)\\[1.3 ex]
\xi(u)\mapsto t_{22}(u)t_{11}(u)^{-1}-
t_{21}(u)t_{11}(u)^{-1}t_{12}(u)t_{11}(u)^{-1}
\end{gather*}
In particular,
\begin{gather*}
\imath(e_0) = t_{21}^{(1)}\qquad\qquad
\imath(f_0) = t_{12}^{(1)}\qquad\qquad
\imath(\xi_0) = t_{22}^{(1)}-t_{11}^{(1)}\\[1 ex]
\imath(\xi_1)
= t_{22}^{(2)}-t_{11}^{(2)}
+\fp\lp(t_{11}^{(1)})^2-t_{22}^{(1)}t_{11}^{(1)}-t_{21}^{(1)}t_{12}^{(1)}\rp
\end{gather*}
which implies that the element $t_1=\xi_1-\fp\xi_0^2/2$ is mapped to
\begin{equation}\label{eq:i t1}
\imath(t_1)
= t_{22}^{(2)}-t_{11}^{(2)}+\frac{\fp}{2}(t_{11}^{(1)}-t_{22}^{(1)})(I+1)-\frac{\fp}{2}\kappa
\end{equation}
where $I=t_{11}^{(1)}+t_{22}^{(1)}$ and $\kappa=t_{12}^{(1)}t_{21}^{(1)}
+t_{21}^{(1)}t_{12}^{(1)}$.

\begin{rem}
The restriction of $\imath$ to $\sl_2\subset\Ysl{2}$ is not the standard embedding
$\jmath:\sl_2\to\gl_2$ given by $e\to E_{12},f\to E_{21},h\to E_{11}-E_{22}$. In fact,
$\left.\imath\right|_{\sl_2}=\theta\circ\jmath$, where $\theta\in\Aut(\gl_2)$ is the
Chevalley involution given by
\begin{equation}\label{eq:theta gl2}
\theta(E_{ij})=E_{\ol{i}\,\ol{j}}
\end{equation}
with $\ol{1}=2$ and $\ol{2}=1$.
\end{rem}

\subsection{The trigonometric Casimir connection of $\gl_2$ \cite{valerio3}}
\label{sec: trig-casimir-sub}

Let $T\subset GL_2$ be the maximal torus consisting of diagonal matrices
and $\tt$ its Lie algebra. The trigonometric Casimir connection of $\gl_2$
is the $Y_\fp\gl_2$--valued connection on $T$ given by
\begin{equation}\label{eq: trig-casimir-gl}
\trigogl{2}=d-\fp
\frac{d(\ve_1-\ve_2)}{e^{\ve_1-\ve_2} - 1}\kappa-d\ve_1\A_1-d \ve_2\A_2
\end{equation}
where
\begin{enumerate}
\item $\{\ve_1, \ve_2\}$ is the basis of $\tt^*$ given by $\ve_i(E_{jj})=\delta
_{ij}$ and $\{d\ve_i\}$ are the corresponding translation--invariant $1$--forms
on $T$.
\item The elements $\A_1,\A_2\in\Ygl{2}$ are given by
\begin{align*}
\A_1 &= 2t_{11}^{(2)} - \fp (t_{11}^{(1)})^2 - \fp t_{11}^{(1)}\\ 			
\A_2 &= 2t_{22}^{(2)} - \fp (t_{22}^{(1)})^2 - \fp t_{22}^{(1)}-\fp\kappa
\end{align*}
\end{enumerate}
Let the symmetric group $\Sym_2$ act on $\Ygl{2}$ by $\sigma(t_{ij}^{(r)})=
t_{\sigma(i),\sigma(j)}^{(r)}$, and regard $\Ysl{2}$ as embedded in $\Ygl{2}$
via \ref{ss:embed}.

\begin{thm}\label{thm: flatness-trig-casimir}\cite[\S 5]{valerio3}
\begin{enumerate}
\item The trigonometric Casimir connection $\trigogl{2}$ is a flat, $\Sym
_2$--equivariant connection on the trivial vector bundle $T\times \Ygl{2}$.
\item The restriction of $\trigogl{2}$ to the maximal torus $H$ of $SL_2$
is the trigonometric Casimir connection $\trigosl{2}$ of $\sl_2$.
\end{enumerate}
\end{thm}

\subsection{Evaluation homomorphism}\label{sec: ev-yangian-gl}

The Yangian $\Ygl{2}$ admits a one--parameter family of algebra
homomorphisms $\ev_a:\Ygl{2}\to U\gl_2[\fp]$ labelled by $a\in\C
[\fp]$, and given by
\[\ev_a(t_{ij}^{(r)})=a^{r-1}E_{ij}\]
Note that this expression continues to make sense, and to define a
homomorphism $\Ygl{2}\to U\gl_{2}[\fp]$ if $a$ is a central element
in $U\gl_{2}[\fp]$.

The evaluation homomorphism of $\Yhgl{2}$ does not restrict to the
one for $\Ysl{2}$ defined in \ref{sec: ev-yangian}. However, the following
holds
\begin{lem}\label{le:shift}
If the evaluation points are related by
\begin{equation}\label{eq:shift}
r=s+\frac{\fp}{2}(I+1)
\end{equation}
the following diagram is commutative
\[\xymatrix{
\Ysl{2} \ar[rd]_{\ev_s} \ar[rr]^{\imath} && \Ygl{2} \ar[ld]^{\theta\circ\ev_r}\\
& U\gl_2[\fp] &}\]
where $\theta$ is the Chevalley involution \ref{eq:theta gl2}.
\end{lem}
\begin{pf}
This clearly holds for the generators $e_0,f_0,\xi_0$ of $\Ysl{2}$, and follows
for the element $t_1$ by comparing
\[\ev_r(\imath(t_1))=(E_{11}-E_{22})(-r+\frac{\fp}{2}(I+1))-\frac{\fp}{2}\kappa\]
where we used \eqref{eq:i t1}, with $\ev_s(t_1)=sh-\fp\kappa/2$.
\end{pf}

\subsection{} 

Now let $(s_1,\ldots,s_k)\in\C[\fp]^k$ and set $r_a=s_a+\frac{\fp}{2}
(I+1)$ for any $1\leq a\leq k$, as in \eqref{eq:shift}. Consider the algebra
homomorphism
\[\ev_{\ul{r}} = \ev_{r_1}\otimes\cdots\otimes \ev_{r_k}\circ \Delta^{(k)}:
\Ygl{2}\to U\gl_2^{\otimes k}[h]\]

\begin{prop}\cite[Prop. 5.6]{valerio3}
\label{prop: trig-casimir-gl-ev}
\begin{enumerate}
\item The image of $\trigogl{2}$ under the evaluation homomorphism 
$\ev_{\ul{r}}$ is the $U\gl_2^{\otimes k}[\fp]$--valued connection given
by
\begin{equation}\label{eq: trig-casimir-gl-ev}
\trigoglev{2}{s}=d-
\fp\frac{d(\ve_1 - \ve_2)}{e^{\ve_1-\ve_2}-1}\Delta^{(k)}(\kappa)
-d\ve_1A_1-\ve_2A_2
\end{equation}
where
\begin{align*}
A_1 &= \sum_a\left(2s_a E_{11}+\fp E_{11}E_{22}\right)^{(a)}
+2\fp\sum_{a<b} E_{12}^{(a)}E_{21}^{(b)}\\
A_2 &= \sum_a\left(2s_a E_{22}+\fp E_{11}E_{22}\right)^{(a)}
-2\fp\sum_{a<b} E_{12}^{(a)}E_{21}^{(b)}-\fp\sum_a\kappa^{(a)}
\end{align*}
\item The restriction of $\trigoglev{2}{s}$ to $H\subset T$ is the image
of the $U\sl_2^{\otimes k}[\fp]$--valued connection $\trigoslev{2}{s}$ of
Proposition \ref {prop: trig-casimir-ev} under the Chevalley involution
$\theta^{\otimes k}$.
\end{enumerate}
\end{prop}
\begin{pf}
(1) By \ref{sec: yangian-gl}, $\Delta(t_{ii}^{(2)})=t_{ii}^{(2)}\otimes 1+1
\otimes t_{ii}^{(2)}+\fp\sum_{i'}t_{ii'}^{(1)}\otimes t_{i'i}^{(1)}$, which
implies that
\[\Delta^{(k)}(t_{ii}^{(2)})=
\sum_a(t_{ii}^{(2)})^{(a)}+
\fp\sum_{a<b}(t_{i\ol{\imath}}^{(1)})^{(a)}(t_{\ol{\imath}i}^{(1)})^{(b)}+
\fp\sum_{a<b}(t_{ii}^{(1)})^{(a)}(t_{ii}^{(1)})^{(b)}\]
where $\ol{1}=2,\ol{2}=1$. Since $\Delta^{(k)}(t_{ii}^{(1)})^2=2\sum_{a
<b}(t_{ii}^{(1)})^{(a)}(t_{ii}^{(1)})^{(b)}+\sum_a ((t_{ii}^{(1)})^{(a)})^2$,
this yields
\[\ev_{\ul{r}}\left(2t_{ii}^{(2)}-\fp(t_{ii}^{(1)})^2-\fp t_{ii}^{(1)}\right)=
\sum_a \left(E_{ii}(2r_a-\fp(E_{ii}+1))\right)^{(a)}+
2\fp\sum_{a<b}(t_{i\ol{\imath}}^{(1)})^{(a)}(t_{\ol{\imath}i}^{(1)})^{(b)}\]
Substituting $r_a=s_a+\half{\fp}(I+1)$ yields the claimed formula for $
A_1=\ev_{\ul{r}}(\A_1)$. The formula for $A_2$ follows from the above,
and the fact that
\[\Delta^{(k)}(\kappa)=\sum_a\kappa^{(a)}+
2\sum_{a<b}\left((t_{12}^{(1)})^{(a)}(t_{21}^{(1)})^{(b)}+
(t_{21}^{(1)})^{(a)}(t_{12}^{(1)})^{(b)}\right)\]

(2) is a direct consequence of Proposition \ref{thm: flatness-trig-casimir}
and Lemma \ref{le:shift}.
\end{pf}

\begin{rems}\hfill
\begin{enumerate}
\item
Since the Chevalley involution is given by conjugating by the matrix 
$\begin{pmatrix}0& i\\i&0\end{pmatrix}\in SL_2$, the application of $\theta
^{\otimes k}$ to the connection $\trigoslev{2}{s}$ yields a connection with
the same monodromy.
\item
As shown in \cite[\S 5.15]{valerio3}, the connection $\trigoglev{2}{k}$
coincides, modulo abelian terms, with the trigonometric dynamical
differential equations for $\gl_2$ considered in \cite{tarasov-varchenko-duality}.
\end{enumerate}
\end{rems}

\section{Affine braid groups}\label{se:braid gps}

\subsection{}

Set
\[\BSL=\pi_1(H_{\reg}/W)\aand\BGL=\pi_1(T_{\reg}/W)\]
The following is well known \cite{cherednik-hecke,macdonald-affine,
etingof-geer} 

\begin{prop}\label{prop: fundamental-gps}\hfill
\begin{enumerate}
\item $\BSL$ is the affine braid group of type $\sfA_1$, and hence
admits the presentation
\[\BSL=\la S_0, S_1|\text{ no relations }\ra\]
\item $\BGL$ can be realised
as the subgroup of the Artin braid group on three strands $B_3$, 
consisting of braids where the first strand is fixed. It has the presentation
\[\BGL=\la \X_1,b|\,b\X_1b\X_1=\X_1b\X_1b\ra\]
\end{enumerate}
\end{prop}

\subsection{}

We describe the generators $S_0,S_1,b,\X_1$ below, together with
the inclusion $\BSL\subset\BGL$ stemming from the $W$--equivariant
embedding $H_{\reg}\subset T_{\reg}$.

Identify to this end the tori $H$ and $T$ with $\nC$ and $(\nC)^2$
respectively, by
\[z\to\begin{pmatrix}z&0\phantom{^{-1}}\\0&z^{-1}\end{pmatrix}
\aand
(z_1,z_2)\to\begin{pmatrix}z_1&0\\0&z_2\end{pmatrix}\]
In terms of these identifications, the inclusion $H\subset T$ is
given by $z\mapsto(z,z^{-1})$. Moreover, $H_{\reg}\subset H$
and $T_{\reg}\subset T$ are identified with $\nC\setminus\{\pm 1\}$
and $Y_2(\nC)$ respectively, where the latter is the configuration
space of two ordered points in $\nC$.

\subsection{}\label{sec: inc-ex-begin}

\begin{figure}
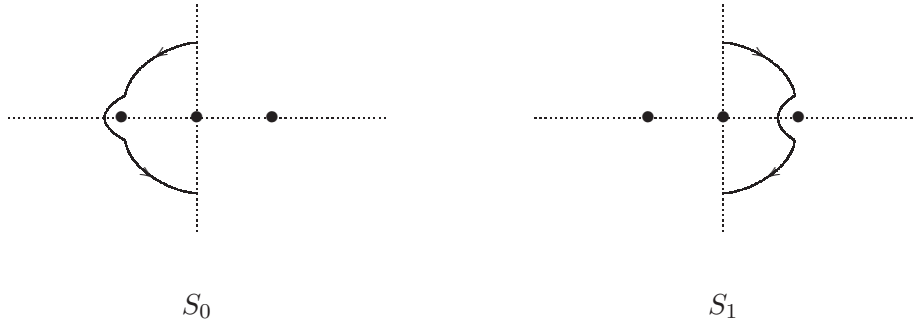

\[
\xy
(0,0)*++{\xy
(-25,0)*{}; (25,0)*{} **\dir{.};
(0,-15)*{}; (0,15)*{} **\dir{.};
(-10,0)*{\bullet}="a"; (10,0)*{\bullet}="b";
(0,10)*{}="c"; (0,-10)*{}="d";
(0,0)*{\bullet};
(-9.5,3)*{}="x"; (-9.5,-3)*{}="y";
"x"; "c" **\crv{(-9,7)&(-3,10)}?(0.5)*\dir{<};
"x"; "y" **\crv{(-15,0)};
"y"; "d" **\crv{(-9,-7)&(-3,-10)}?(0.5)*\dir{>};
(0,-25)*{S_0};
\endxy}="P";
(70,0)*++{\xy
(-25,0)*{}; (25,0)*{} **\dir{.};
(0,-15)*{}; (0,15)*{} **\dir{.};
(-10,0)*{\bullet}="a"; (10,0)*{\bullet}="b";
(0,10)*{}="c"; (0,-10)*{}="d";
(0,0)*{\bullet};
(9.5,3)*{}="x"; (9.5,-3)*{}="y";
"x"; "c" **\crv{(9,7)&(3,10)}?(0.5)*\dir{<};
"x"; "y" **\crv{(5,0)};
"y"; "d" **\crv{(9,-7)&(3,-10)}?(0.5)*\dir{>};
(0,-25)*{S_1};
\endxy}="Q";
\endxy
\]
\caption{Generators of $\pi_1(H_{\reg}/W)$}\label{fig: reg-torus}
\end{figure}

The generators $S_0, S_1$ of $\BSL$ may be described as follows
\cite{nguyen, van-der-lek-2, van-der-lek}. Identify the Lie algebra
$\h$ of $H$ with $\C$ by mapping $h$ to $1$. The exponential map
$\exp\lp 2\pi\iota-\rp:\h\to H$ maps $\h^{\areg}$ to $H_{\reg}$, where
\[\h^{\areg}=\h\setminus\bigcup_{n\in \Z}\{\alpha=n\} \cong \C \setminus\frac{1}{2}\,\Z\]
The affine Weyl group $W_{\aff}$ of type $\sfA_1$ is generated by the
affine (real) reflections $s_0, s_1$ through the points $u=1/2$ and $u=0$
respectively. $W_{\aff}$ is isomorphic to $\Z_2\ltimes \Z$, with the generator
$s_1$ of $\Z_2$ acting on $\h$ as the reflection $u\to -u$ and the generator
$\tau=s_0s_1$ of $\Z$ as the translation $u\to u+1$. Thus we have the
identification 
\begin{equation}\label{eq:identification}
\exp(2\pi\iota-) : \h^{\areg}/W_{\aff} \cong H_{\reg}/W
\end{equation}
Fix now a base point (say $u=1/4$) in $\h^{\areg}$ lying in the interval
$(0,1/2)$. Then, the generators $S_i$ are represented by the loops in
$\h^{\areg}/W_{\aff}$ given in Figure \ref{fig: aff-reg}.
\begin{figure}
\[
\xy 0;/r.14pc/:
(0,0)*++{\xy
(-40,0)*{\bullet}; (-20,0)*{\bullet}; (0,0)*{\bullet}="a";
(20,0)*{\bullet}="b"; (40,0)*{\bullet};
(5,0)*{\cdot}="x"; (35,0)*{\cdot}="y";
(13,0)*{}="c"; (27,0)*{}="d";
(-50,0)*{}; (50,0)*{} **\dir{.};
(0,-20)*{}; (0,20)*{} **\dir{.};
"x"; "c" **\crv{}?(0.6)*\dir{>};
"d"; "y" **\crv{}?(0.6)*\dir{>};
(14,0)*{}; (26,0)*{} **\crv{(18,-6)&(22,-6)};
(0,-25)*{S_0};
\endxy}="P";
(140,0)*++{
\xy
(-40,0)*{\bullet}; (-20,0)*{\bullet}="b"; (0,0)*{\bullet}="a";
(20,0)*{\bullet}; (40,0)*{\bullet};
(15,0)*{\cdot}="x"; (-15,0)*{\cdot}="y";
(7,0)*{}="c"; (-7,0)*{}="d";
(-50,0)*{}; (50,0)*{} **\dir{.};
(0,-20)*{}; (0,20)*{} **\dir{.};
"x"; "c" **\crv{}?(0.6)*\dir{>};
"d"; "y" **\crv{}?(0.6)*\dir{>};
(6,0)*{}; (-6,0)*{} **\crv{(2,6)&(-2,6)};
(0,-25)*{S_1};
\endxy}="Q";
\endxy
\]
\caption{Generators of $\pi_1\lp\h^{\areg}/W_{\aff}\rp$}\label{fig: aff-reg}
\end{figure}
These correspond, via the identification \eqref{eq:identification} to the loops
in $H_{\reg}$ shown in Figure ~\ref{fig: reg-torus}.

\subsection{}\label{sec: pi-2}

Turning to the fundamental group $\BGL$, it is conventional
to pick its base point as the configuration $(1,2)$ in $\nC$.
The generators $\X_1,b$ are then represented by the braids
in Figure ~\ref{fig: pi-2}. Set $\X_2=b\X_1b$. Then, the
defining relation of $\BGL$ can be written as $\X_1\X_2=\X_2
\X_1$.

\begin{figure}
\[
\xy
(0,0)*++{\xy
(0,-3)*{0}; (10,-3)*{1}; (20,-3)*{2}; 
(0,23)*{0}; (10,23)*{1}; (20,23)*{2}; 
(0,0)*{}="o1"; (10,0)*{}="a"; (20,0)*{}="b"; 
(0,20)*{}="o2"; (10,20)*{}="g"; (20,20)*{}="h"; 
"o1"; "o2" **\dir{-};
"h"; "a" **\dir{-}\POS?(0.5)*{\hole}="x";
"x"; "b" **\dir{-};
"x"; "g" **\dir{-};
(-5,0)*{}; (25,0)*{} **\dir{-};
(-5,20)*{}; (25,20)*{} **\dir{-};
(10,-10)*{b};
\endxy}="P";
(50,0)*++{\xy
(0,-3)*{0}; (10,-3)*{1}; (20,-3)*{2};
(0,23)*{0}; (10,23)*{1}; (20,23)*{2}; 
(0,0)*{}="a"; (10,0)*{}="b"; (20,0)*{}="c"; 
(0,20)*{}="g"; (10,20)*{}="h"; (20,20)*{}="i"; 
(-3,0)*{}; (23,0)*{} **\dir{-};
(-3,20)*{}; (23,20)*{} **\dir{-};
"c"; "i" **\dir{-};
"g"; (0,10)*{} **\dir{-}\POS?(0.8)*{\hole}="x";
"h"; "x" **\crv{(10,16)};
"x"; "b" **\crv{(-3,11)&(-3,9)&(0,8)&(10,4)};
"a"; (0,7)*{} **\dir{-};
(10,-10)*{\X_1};
\endxy}="Q";
\endxy
\]
\caption{Generators of $\BGL$}\label{fig: pi-2}
\end{figure}

\subsection{}\label{sec: inc-ex-end}

To relate $\BSL$ and $\BGL$, think of elements of $\BSL$ as braids
with 3 strands, with endpoints $-i,0,i$, and the strand at zero remaining
fixed. Choosing a path from the base point $(-i,0,i)$ to $(0,1,2)$ which
first braids the first two points to $(0,i/2,i)$ while keeping the third fixed
and then scales the configuration to $(0,1,2)$ yields an embedding $\BSL\to\BGL$ given by
\begin{equation}\label{eq: inclusion}
\begin{array}{ccl}
S_1 & \mapsto & b \\
S_0 & \mapsto & \X_1 b\X_1^{-1}
\end{array}
\end{equation}

Let $\LL=S_0S_1$ be the element of $\BSL$ corresponding to the
generator of the coroot lattice of $SL_2$. Then, the inclusion above
yields
\[\LL \mapsto \X_1b\X_1^{-1}b= b^{-1}(b\X_1b) \X_1^{-1} b= b^{-1}(\X_2\X^{-1})b\]
Thus, if we consider the set $\{b,\LL_1,\LL_2\}$ of generators of $\BGL$
obtained by conjugation with $b^{-1}$
\[\LL_1=b^{-1}\X_2b\aand\LL_2=b^{-1}\X_1b\]
then, the image of the element $\LL$ of $\BSL$ is given by
\begin{equation}\label{eq: inclusion-lattice}
\LL \mapsto \LL_1\LL_2^{-1}
\end{equation}

\subsection{}

We now relate the monodromy representations of the trigonometric
Casimir connections $\trigoslev{2}{s}$ and $\trigoglev{2}{s}$. Note first that
the actions of $W \cong \Z/2\Z$ on the fibers of the corresponding
vector bundles are different. This difference arises from the fact that in the
case of $SL_2(\C)$,  the following element is used to construct a group
homomorphism $\wt{W} \to SL_2(\C)$ (see the discussion preceding
\cite[Corollary 3.6]{valerio4}).
\[\sigma \mapsto \left(\begin{array}{rr} 0 & 1 \\ -1 & 0 \end{array}\right)\]
Let $\pisl$ and $\pigl$ denote the  representations of $\BSL$ and $\BGL$
obtained from the monodromy of the connections $\trigoslev{2}{s}$ and
$\trigoglev{2}{s}$ respectively. Then we have
\begin{equation}\label{thm: monodromy-sl-gl}
\begin{aligned}
\pisl(b)&=\pigl(b)(-1)^{E_{11}}\\
\pisl(\LL)&=\pigl(\LL_2\LL_1^{-1})
\end{aligned}
\end{equation}

\section{The trigonometric KZ equations}\label{se:trigo KZ}

\subsection{}\label{ss: trig-KZ}

Fix $k\geq 2$, let $r\in\gl_k^{\otimes 2}$ be the Drinfeld $r$--matrix of $\gl_k$,
\[r=\frac{1}{2}\sum_{a=1}^k\El{aa}\otimes\El{aa}
+\sum_{1\leq a<b\leq k} \El{ab}\otimes \El{ba}\]
and let $r(u)=\ds \frac{re^u + r_{21}}{e^u-1}$ be the corresponding trigonometric
$r$--matrix. Fix $n\geq 1$, let $V$ be a $\gl_k$--module and $\IV^{\otimes n}$
the trivial vector bundle over $\C^n$ with fibre $V^{\otimes n}$. The symmetric
group $\Sym_n$ acts both on the base and fibre of $\IV^{\otimes n}$. The trigonometric
KZ connection is the flat, $\Sym_n$--equivariant connection on $\IV^{\otimes n}$
given by
\[\nabla\KKZ=d-2\fp\left(\sum_{i<j}r_{ij}(u_i-u_j)d(u_i-u_j)+\sum_i s^{(i)}\,du_i\right)\]
where $s\in\gl_k$ is a fixed diagonal matrix and $s^{(i)}=1^{\otimes (i-1)}
\otimes s\otimes 1^{\otimes (n-i)}$. Since $r(u)=\Omega/(e^u-1)+r$, where
$\Omega=r+r_{21}$, this connection may equivalently be written as
\begin{equation}\label{eq:trig-KZ-eg}
\nabla\KKZ=d-2\fp\left(
\sum_{i<j}\frac{d(u_i-u_j)}{e^{u_i-u_j}-1}\,\Omega_{ij}
+\sum_i du_i\,X_i\right)
\end{equation}
where $X_i=s^{(i)}+\sum_{j>i}r_{ij}-\sum_{j<i}r_{ji}$.

\subsection{}\label{ss:Pi_n}

The connection $\nabla\KKZ$ is invariant under the group $\Z^n$ acting
trivially on the fibres of $\IV^{\otimes n}$ and by translations by the lattice $2
\pi\imath\Z^n$ on the base. It therefore descends to a flat connection on the
complement $\XX_n$, in the quotient $\C^n/\Sym_n\ltimes\Z^n$ of the images
of the affine hyperplanes $\{u_i-u_j=2\pi\imath m\}_{i\neq j,m\in\Z}$. The latter
may be thought of as the configuration space of $n$ points in $\nC$, or equivalently
the set of regular elements in the maximal torus of diagonal matrices in $GL_n
(\C)$. The following gives a presentation of the fundamental group $\Pi_n=B_
{GL_n}$ of this space

\begin{prop}\cite{birman-braids}
$\Pi_n$ is generated by elements $\{b_i\}_{1\leq i\leq n-1}$ and
$\{X_j \}_{1\leq j\leq n}$, subject to the relations
\begin{alignat*}{3}
b_ib_{i'}		&=b_{i'}b_i		\\
b_ib_{i+1}b_i	&=b_{i+1}b_ib_{i+1}	\\
b_i\X_ib_i		&=\X_{i+1}		\\
\X_j\X_k		&=\X_k\X_j		
\end{alignat*}
for any $1\leq i,i'\leq n-1$ such that $|i-i'|\geq 2$, and $1\leq j,k\leq n$.
\end{prop}

The generators $b_i,\X_j$ may be described as follows. Let $z_i=e^{u_i}$, $i=1,
\ldots,n$, be the standard coordinates on $\XX_n$ and, for definiteness, choose
$\ul{z}_0=(1,\ldots,n)$ as basepoint. Then, $\X_j$ and $b_i$ are, respectively,
the loops
\begin{align*}
t&\mapsto (1,\ldots,j-1,e^{2\pi\imath t}j,j+1,\ldots,n)\\
t&\mapsto (1,\ldots,i-1,i+1/2(1-e^{\pi\imath t}),i+1/2(1+e^{\pi\imath t}),i+2,\ldots,n)
\end{align*}
where $t\in[0,1]$.

\section{The dual pair $(\gl_k,\gl_2)$ and trigonometric connections}\label{sec: howe}

\subsection{}\label{sec: dual-pair}

Let $\M_{k,2}$ be the vector space of complex, $k\times 2$ matrices
and
\[\cfun{k}{2}=\C[x_{aj}]_{\substack{1\leq a\leq k\\[.2 ex] 1\leq j\leq 2}}\]
its algebra of regular functions. The group $GL_k\times GL_2$ acts
on $\cfun{k}{2}$ by
\[(g_k,g_2)p(X)=p(g_k^t\,X\,g_2)\]
where $X\in\M_{k,2}$ and $g_p\in GL_p$. Note that
\begin{equation}\label{eq:ids}
\C[x_{a1}]\otimes\C[x_{a2}]\cong\cfun{k}{1}^{\otimes 2}\cong
\cfun{k}{2}\cong
\cfun{1}{2}^{\otimes k}\cong\C[x_{1j}]\otimes\cdots\otimes\C[x_{kj}]
\end{equation}
where the first two are isomorphisms of $GL_k$--modules, and the
last two of $GL_2$--modules.

The action of $GL_k\times GL_2$ preserves the \fd homogeneous
components of $\cfun{k}{2}$ and therefore gives rise to an action of
$\gl_k\oplus\gl_2$ on these. To distinguish between the elements of
these Lie algebras, we denote by $X^{(p)}$ the elements of $\gl_p$.
Then, the action is given by mapping the elementary matrices $\el{ab}
{k},\el{ij}{2}$ to
\begin{align}\label{eq: joint-action}
\el{ab}{k} &\mapsto \sum_{j=1}^2 x_{aj}\partial_{bj}
& \el{ij}{2} &\mapsto \sum_{a=1}^k x_{ai}\partial_{aj}
\end{align}

\subsection{}

\begin{lem}\label{le:internal}
The following holds on $\cfun{1}{2}$
\[\kappa^{(2)}=I^{(2)}+2E_{11}^{(2)}E_{22}^{(2)}\]
\end{lem}
\begin{pf}
On $\cfun{1}{2}\cong\C[x_j]_{j=1,2}$, $\kappa=E_{12}E_{21}+E_{21}E_{12}$ acts as
\[x_1\partial_2x_2\partial_1+x_2\partial_1x_1\partial_2=
x_1\partial_1x_2\partial_2+x_1\partial_1+
x_2\partial_2x_1\partial_1+x_2\partial_2=
2E_{11}E_{22}+I\]
\end{pf}

\subsection{Duality}

The identities below relate the coefficients of the KZ connection of $\gl_k$ and
those of the Casimir connection of $\gl_n$ (here, $n=2$). They were discovered
in \cite{valerio4}. Let $r=r^{(k)}\in\gl_k^{\otimes 2}$ be the $r$--matrix defined in
\S \ref{ss: trig-KZ} and $\Omega^{(k)}=r+r_{21}$. 

\begin{prop}\label{pr:transfer}
The following identities hold on $\cfun{k}{1}^{\otimes 2}\cong\cfun{k}{2}\cong\cfun{1}{2}
^{\otimes k}$
\begin{gather*}
(E_{aa}^{(k)})^{(i)}				=(E_{ii}^{(2)})^{(a)}\\[1.2 ex]
r^{(k)}						=\sum_{a<b}(E_{12}^{(2)})^{(a)}(E_{21}^{(2)})^{(b)}+
								\half{1}\sum_a(E_{11}^{(2)}E_{22}^{(2)})^{(a)}\\[1.2 ex]
2\Omega^{(k)}					=\Delta^{(k)}(\kappa^{(2)}-I^{(2)})\\[1.2 ex]
s^{(i)	}						=\sum_a s_a(E_{ii}^{(2)})^{(a)}\\[1.2 ex]
\end{gather*}
\end{prop}
\begin{pf}
(1) The identity follows from the fact that both sides record the homogeneity
degree \wrt the variable $x_{ai}$.

(2) By \eqref{eq: joint-action}, the action of $\ol{r}=r-\half{1}\sum_a E_{aa}^{(k)}\otimes
E_{aa}^{(k)}$ is given by
\[\sum_{a<b}x_{a1}\partial_{b1}x_{b2}\partial_{a2}=
\sum_{a<b}x_{a1}\partial_{a2}x_{b2}\partial_{b1}=
\sum_{a<b}(E_{12}^{(2)})^{(a)}(E_{21}^{(2)})^{(b)}\]
and, by (1), $r-\ol{r}$ acts by $\half{1}\sum_a (E_{11}^{(2)}E_{22}^{(2)})^{(a)}$.

(3) By (2),
\[\begin{split}
\Omega^{(k)}
&=
r+r^{21}\\
&=
\sum_{a<b}\left((E_{12}^{(2)})^{(a)}(E_{21}^{(2)})^{(b)}+
(E_{21}^{(2)})^{(a)}(E_{12}^{(2)})^{(b)}\right)+\sum_a (\el{11}{2}\el{22}{2})^{(a)}\\
&=
\sum_{a\neq b}(E_{12}^{(2)})^{(a)}(E_{21}^{(2)})^{(b)}+\half{1}\sum_a (\kappa^{(2)}-I^{(2)})^{(a)}
\end{split}\]
where we used Lemma \ref{le:internal}. On the other hand, since $\kappa^{(2)}=
\el{12}{2}\el{21}{2}+\el{21}{2}\el{12}{2}$,
\[\begin{split}
\Delta^{(k)}(\kappa^{(2)})
&=
\sum_{a,b}(\el{12}{2})^{(a)}(\el{21}{2})^{(b)}+(\el{21}{2})^{(b)}(\el{12}{2})^{(a)}\\
&=
2\sum_{a\neq b}(\el{12}{2})^{(a)}(\el{21}{2})^{(b)}+
\sum_a(\kappa^{(2)})^{(a)}
\end{split}\]

(4) follows from (1) since
\[s^{(i)}=\sum_a s_a(E_{aa}^{(k)})^{(i)}=\sum_a s_a (E_{ii}^{(2)})^{(a)}\]
\end{pf}

\subsection{}

The following is a direct consequence of Proposition \ref{pr:transfer}
(see also \cite{tarasov-varchenko-duality}).

\begin{prop}\label{pr:duality}
Under the identification
\[\cfun{k}{1}^{\otimes 2}\cong\cfun{k}{2}\cong\cfun{1}{2}^{\otimes k}\]
the trigonometric KZ connection for $\gl_k$ with values in $\cfun{k}{1}
^{\otimes 2}$ corresponding to a diagonal matrix $s=\sum_a s_a\el{aa}
{k}$, coincides with the sum of
\begin{enumerate}
\item the trigonometric Casimir connection for $\gl_2$ with values in the
tensor product of evaluation modules
\[\cfun{1}{2}(r_1)\otimes\cdots\otimes\cfun{1}{2}(r_k)\]
where $r_a=\fp\left(s_a+\half{I^{(2)}+1}\right)$, and 
\item the closed one--form with values in $Z(\gl_2)^{\otimes k}$ given by
\[\A=\fp\left(
\frac{d(\epsilon_1-\epsilon_2)}{e^{\epsilon_1-\epsilon_2}-1}-
d\epsilon_2\right)\,\Delta^{(k)}(I^{(2)})\]
\end{enumerate}
\end{prop}
\begin{pf}
Using the form \eqref{eq:trig-KZ-eg}, it follows from Proposition \ref{pr:transfer}
that the trigonometric KZ connection for $\gl_k$ may be rewritten as the $U\gl_
2^{\otimes k}[\fp]$--valued connection
\[d-\fp\left(
\frac{d(\epsilon_1-\epsilon_2)}{e^{\epsilon_1-\epsilon_2}-1}\,\Delta^{(k)}(\kappa^{(2)})
+d\epsilon_1\,X_1+d\epsilon_2\,X_2\right)+
\fp\frac{d(\epsilon_1-\epsilon_2)}{e^{\epsilon_1-\epsilon_2}-1}\,\Delta^{(k)}(I^{(2)})\]
where
\begin{align*}
X_1 &= 2(s^{(1)}+r)=\sum_a(2s_aE_{11}+\el{11}{2}\el{22}{2})^{(a)}+
2\sum_{a<b}(\el{12}{2})^{(a)}(\el{21}{2})^{(b)}\\
X_2 
&= 2(s^{(2)}-r)
=
\sum_a(2s_a\el{11}{2}-\el{11}{2}\el{22}{2})^{(a)}
-2\sum_{a<b}(\el{12}{2})^{(a)}(\el{21}{2})^{(b)}\\
&=
\sum_a(2s_a\el{11}{2}+\el{11}{2}\el{22}{2})^{(a)}
-2\sum_{a<b}(\el{12}{2})^{(a)}(\el{21}{2})^{(b)}
-\sum_a (\kappa^{(2)}-I^{(2)})^{(a)}
\end{align*}			
where we used Lemma \ref{le:internal}. The result now follows from Proposition
\ref{prop: trig-casimir-gl-ev}.
\end{pf}

\subsection{}

\begin{cor}\label{cor: monodromy-comparison}
Let
\[\pi_{\KZ, \ul{s}},\pi_{C,s}:\BGL\to GL\lp\cfun{k}{2}[[\fp]]\rp\]
be the monodromy representations of the trigonometric KZ
connection for $\gl_k$ and Casimir connection for $\gl_2$
corresponding to the diagonal matrix $s=\sum_a s_a\el{aa}{k}$
and evaluation points $r_a=\fp(s_a+\half{I^{(2)}+1})$ respectively.
Then,
\[\begin{array}{ccl}
\pi_{\KZ,\ul{s}}(b) 	&=& \pi_{C,\ul{s}}(b)e^{-\pi\iota\fp \Delta^{(k)}(I^{(2)})}\\[1.1 ex]
\pi_{\KZ,\ul{s}}(\X_1) &= & \pi_{C,\ul{s}}(\X_1)e^{2\pi\iota\fp \Delta^{(k)}(I^{(2)})}\\[1.1 ex]
\pi_{\KZ,\ul{s}}(\X_2) &= & \pi_{C,\ul{s}}(\X_2)
\end{array}\]
\end{cor}
\begin{pf}
In terms of the coordinates $z_1=e^{\ve_1}$, $z_2=e^{\ve_2}$, the 1--form
$\A$ of Proposition \ref{pr:duality} is equal to
\[\fp\left(\frac{d(z_1-z_2)}{z_1-z_2}-\frac{dz_1}{z_1}-\frac{dz_2}{z_2}\right)\,
\Delta^{(k)}(I^{(2)})\]
A fundamental solution of the corresponding connection is given by
\[\left(z_1z_2/(z_2-z_1)\right)^{\fp\Delta^{(k)}(I^{(2)})}=
\exp\left(\fp\Delta^{(k)}(I^{(2)})(\log z_1+\log z_2 - \log (z_2-z_1))\right)\]
where $\log$ is the standard determination of the logarithm, and has
monodromy along the generators $b,\X_1,\X_2$ described in Section
\ref{ss:Pi_n} given by
\[b\mapsto e^{-\pi\iota\fp \Delta^{(k)}(I^{(2)})}\qquad
\X_1\mapsto e^{2\pi\iota\fp \Delta^{(k)}(I^{(2)})}\qquad
\X_2\mapsto 1\]
\end{pf}

\section{Monodromy of the trigonometric KZ equations}
\label{sec: etingof-geer}

In this section we recall the main theorem of \cite{etingof-geer}
which computes the monodromy of the trigonometric KZ equations.

\subsection{The quantum group $\Uhgl{p}$}\label{sec: qg-gl}

The Drinfeld--Jimbo quantum group $\Uhgl{p}$ is defined as a unital 
associative $\C[[\hbar]]$--algebra, topologically generated by elements
$\{E_j, F_j\}_{1\leq j\leq p-1}$ and $\{D_i\}_{1\leq i\leq p}$ subject to
the relations (where $q^2=e^{\hbar}$)
\begin{itemize}
\item[(QG1)] $[D_i, D_{i'}]=0$ for any $i,i'$.
\item[(QG2)] For each $i, j$, $1\leq i\leq p$ and $1\leq j\leq p-1$ we 
have
\begin{align*}
[D_i,E_j] &= \lp\delta_{ij}-\delta_{i,j+1}\rp E_j
& [D_i, F_j] &= \lp \delta_{i,j+1} - \delta_{ij}\rp F_j
\end{align*}
\item[(QG3)] For each $j,j'\in \{1,\cdots, p-1\}$ we have
\[
[E_j,F_{j'}]=\delta_{j,j'} \frac{q^{H_j} - q^{-H_j}}{q - q^{-1}}
\]
\item[(QG4)] For each $j\neq j' \in \{1,\cdots, p-1\}$ we have:
\begin{gather*}
\sum_{t=0}^{1-a_{jj'}} (-1)^t \qbin{1-a_{jj'}}{t}{q} E_i^{1-a_{jj'}-t} 
E_j E_i^t=0\\
\sum_{t=0}^{1-a_{jj'}} (-1)^t \qbin{1-a_{jj'}}{t}{q} F_i^{1-a_{jj'}-t} 
F_j F_i^t=0
\end{gather*}
\end{itemize}
where $H_i=D_i - D_{i+1}$. We have used the standard notations of
 the Gaussian integers.
\begin{gather*}
[n]_q=\frac{q^n-q^{-n}}{q-q^{-1}} \aand [n]_q!=[n]_q[n-1]_q\cdots
 [1]_q \\
\qbin{n}{m}{q}=\frac{[n]_q!}{[m]_q![n-m]_q!}
\end{gather*}
and $a_{jj'}$ are entries of the Cartan matrix of type $A_{p-1}$:
\[
a_{jj'}=2- \delta_{[j-j'|=1}
\]

$\Uhgl{p}$ is a topological Hopf algebra with coproduct and counit 
given by:
\begin{equation}\label{eq: hopf-Uhgl}
\begin{array}{ccl}
\Delta(D_i) &=& D_i\otimes 1 + 1\otimes D_i\\
\Delta(E_j) &=& E_j\otimes q^{H_j} + 1\otimes E_j\\
\Delta(F_j) &=& F_j\otimes 1 + q^{-H_j}\otimes F_j
\end{array}
\end{equation}
and
\begin{gather}\label{eq: counit-Uhgl}
\ve(E_j) =\ve(F_j)=\ve(D_i)=0
\end{gather}

Let $I_p=D_1 + \cdots + D_p \in \Uhgl{p}$. It is clear from the definition
 above that $I_p$ is a central element of $\Uhgl{p}$, and the coproduct 
 on $I_p$ is given by $\Delta(I_p)=I_p\otimes 1 + 1\otimes I_p$. We 
 have the following isomorphism of Hopf algebras:
\[
\Uhgl{p}=\Uhsl{p}\otimes \C[I_p][[\hbar]]
\]

Moreover $\Uhgl{p}$ has a \qt structure. Let $\RR$ be 
the $R$--matrix of $\Uhgl{p}$. Recall that the Drinfeld element $u$ is 
defined by:
\begin{equation}\label{eq: drinfeld-element}
u=\lp m\circ (S\otimes 1) \rp (\RR_{21})
\end{equation}
The following theorem is proved in \cite{drinfeld-almost}.

\begin{thm}\label{thm: drinfeld-element}
The square of the antipode is an inner automorphism given by:
\[
S^2(x)=uxu^{-1}
\]
\end{thm}

\begin{rem}\label{rem: R-matrix-comparison}
In this note $\RR$ denotes the $R$--matrix of $\Uhgl{p}$, which differs
 from the $R$--matrix of $\Uhsl{p}$ (the one used in \cite{valerio4}) by:
\begin{equation}\label{eq: R-matrix-comparison}
\RR_{\gl_p}=q^{\frac{I_p\otimes I_p}{p}} \RR_{\Lsl_p}
\end{equation}
\end{rem}

\subsection{Monodromy of the trigonometric KZ equations \cite{etingof-geer}}
\label{sec: monodromy-trig-KZ}

Let $V$ be a $\gl_k$--module on which $I$ acts semisimply, $n\geq 1$
and consider the monodromy of the trigonometric KZ equations defined
in Section \ref{se:trigo KZ} on $V^{\otimes n}$.

Let $\V$ be a \fd $\Uhgl{k}$--module satisfying $\V/\hbar\V \cong V$
and such that $I$ acts semisimply and with eigenvalues in $\C$. Define
\begin{equation}\label{eq: eg-element}
T=(S\otimes\Id)(\RR_{21}) \aand 
C=m_{01}\lp T_{0n}\cdots T_{01}\rp=m_{01}\lp(1\otimes\Delta^{(n)})T\rp
\end{equation}
The following is the main result of \cite{etingof-geer}. It relies upon the fact
that the \EK quantization of $\gl_k$ corresponding to the $r$--matrix given
in \ref{ss: trig-KZ} coincides with the Drinfeld--Jimbo quantum group $\Uhgl
{k}$, which is proved in \cite{etingof-kazhdan-quantization-6}. 

\begin{thm}\label{thm: etingof-geer}
Let $\hbar=4\pi\iota h$. Then, the monodromy representation $\pi : \Pi_n \to
GL\lp V^{\otimes n}[[\hbar]]\rp$ corresponding to \eqref{eq:trig-KZ-eg} is
equivalent to  the following representation of $\Pi_n$ on $\V^{\otimes n}$
\begin{align*}
b_i &\mapsto (i\ i+1)\RR_{i,i+1} \\
\X_1 & \mapsto (q^{2s}u^{-1})^{(1)}\,C
\end{align*}
\end{thm}
The statement of the above theorem differs slightly from the one given
in \cite{etingof-geer}, due to minor computational errors in \cite{etingof-geer}.
For reader's convenience, we reproduce the proof of this theorem in
Appendix \ref{app-sec: etingof-geer}.

\subsection{}
\label{sec: cor-eg}

\begin{cor}\label{cor: cor-eg}
The monodromy of the trigonometric KZ connection \eqref{eq:trig-KZ-eg}
is equivalent to the action of $\Pi_n$ on $\V^{\otimes n}$ given by
\begin{align*}
\rho_{s}(b_i)	&=(i\,i+1)\RR_{i\,i+1}\\
\rho_{s}(\X_j)	
&=
\RR_{j\,j-1}\cdots\RR_{j\,1}\,
(q^{2\fs})^{(j)}\,
\RR_{n\,j}^{-1}\cdots\RR_{j+1\,j}^{-1}\\
&=
\Delta^{(j-1)}\otimes\Id(\RR_{21})\cdot (q^{2\fs})^{(j)}\cdot
1^{\otimes (j-1)}\otimes\Id\otimes\Delta^{(n-j)}(\RR_{21}^{-1})
\end{align*}
\end{cor}
\begin{pf}
We need only check the assignment for $\X_1,\ldots,\X_n$. Write $\RR=
\alpha_i\otimes\beta^i$, where the sum over $i$ is implicit. By \eqref{eq:
eg-element},
\[\begin{split}
(u^{-1})^{(1)}C
&=(u^{-1})^{(1)}m_{01}(T_{0n}\cdots T_{01})\\
&=
u^{-1}S(\beta^{i_n})\cdots S(\beta^{i_1})\alpha_{i_1}\otimes\alpha_{i_2}\otimes\cdots\otimes\alpha_{i_n}\\
&=
u^{-1}S(\beta^{i_n})\cdots S(\beta^{i_2})u\otimes\alpha_{i_2}\otimes\cdots\otimes\alpha_{i_n}\\
&=
S^{-1}(\beta^{i_2}\cdots\beta^{i_n})\otimes\alpha_{i_2}\otimes\cdots\otimes\alpha_{i_n}\\
&=
S^{-1}\otimes\Delta^{(n-1)}(\RR_{21})\\
&=
\Id\otimes\Delta^{(n-1)}(\RR_{21}^{-1})
\end{split}\]
where we used $u=S(\beta^{i_1})\alpha_{i_1}$, $\Ad(u)=S^2$, the cabling
identity $\Delta^{(n-1)}\otimes\Id(\RR)=\RR_{1\,n}\RR_{2\,n}\cdots\RR_{n-1\,n}$
which implies that
\[\Id\otimes\Delta^{(n-1)}(\RR_{21})=
(1\,2\cdots n)\circ\Delta^{(n-1)}\otimes\Id(\RR)=
\RR_{2\,1}\RR_{3\,1}\cdots\RR_{n\,1}\]
and the fact that $(S^{-1}\otimes 1)(\RR_{21})=\RR_{21}^{-1}$. Thus,
\[\rho_{s}(\X_1)=
(q^{2\fs})^{(1)}\,\Id\otimes\Delta^{(n-1)}(\RR_{21}^{-1})=
(q^{2\fs})^{(1)}\,
\RR_{n\,1}^{-1}\cdots\RR_{2\,1}^{-1}\]
The formula for $\X_j$ follows from this by an easy induction using $\X_{j+1}
=b_j\X_jb_j$.
\end{pf}

\subsection{}\label{ss:EG n=2} 

We shall mostly be interested in the case $n=2$. In this case,
the above formulae read
\[\rho_{s}(b)=(1\,2)\RR\qquad
\rho_{s}(\X_1)=(q^{2\fs})^{(1)}\RR_{21}^{-1}\qquad
\rho_{s}(\X_2)=\RR_{21}(q^{2\fs})^{(2)}\]
which, in terms of the generators $b,\LL_1=b^{-1}\X_2b,\LL_2
=b^{-1} \X_1b$ of $B_{GL_2}$ defined in \ref{sec: pi-2}, yields
\begin{equation*}\label{eq:EG for BGL}
\rho_{s}(b)		 =(1\,2)\RR\qquad
\rho_{s}(\LL_1) 	 =(q^{2\fs})^{(1)}\RR\qquad
\rho_{s}(\LL_2)	=\RR^{-1}(q^{2\fs})^{(2)}
\end{equation*}

\section{Quantum loop algebras}\label{sec: qla}

In this section we review the definitions of the quantum loop algebras
$\hloopgl{2}$ and $\hloopsl{2}$ following \cite{ding-frenkel} and \cite
{chari-pressley-qaffine} respectively.

\subsection{The quantum loop algebra $\hloopgl{2}$ \cite{ding-frenkel}}
\label{sec: defn-qla}

$\hloopgl{2}$ is a unital, associative, complete $\C[[\hbar]]$--algebra
topologically generated by elements $\{E_r,F_r,D_{1,r},D_{2,r},\}_{r
\in\Z}$. To state its defining relations, consider the formal series
\begin{gather*}
E(z)=\sum_{r\in \Z} E_rz^{-r}\qquad F(z)=\sum_{r\in \Z} F_rz^{-r}\\
\Theta_j^{\pm}(z)=q^{\pm D_{j,0}}\exp\lp \pm (q-q^{-1}) \sum_{r\geq 1}
 D_{j,\pm r} z^{\mp r}\rp
\end{gather*}
Then the relations can be written as\footnote{The presentation above
differs from the one given in \cite{ding-frenkel} by the interchange $\Theta
_1(z)^{\pm}\leftrightarrow\Theta_2(z)^\pm$. The present convention
makes the formulae for the inclusion $\Uhgl{2}\hookrightarrow\hloopgl
{2}$ and evaluation $\hloopgl{2}\to\Uhgl{2}$ more natural.}
\begin{itemize}
\item[(QL1)] The elements $\{D_{j,r}\}_{j=1,2;r\in\Z}$ commute. 
\item[(QL2)] Let $\theta_m(\zeta)=\frac{q^m\zeta-1}{\zeta-q^m}$, then

\begin{alignat*}{3}
\Theta_1^{\pm}(z)E(w)\Theta_1^{\pm}(z)^{-1}&=\theta_{1}(qz/w)E(w)&\qquad
\Theta_2^{\pm}(z)E(w)\Theta_2^{\pm}(z)^{-1}&=\theta_{-1}(q^{-1}z/w)E(w)\\[1.1 ex]
\Theta_1^{\pm}(z)^{-1}F(w)\Theta_1^{\pm}(z)&=\theta_{1}(qz/w)F(w)&\qquad
\Theta_2^{\pm}(z)^{-1}F(w)\Theta_2^{\pm}(z)&=\theta_{-1}(q^{-1}z/w)F(w)
\end{alignat*}
\item[(QL3)]
\begin{align*}
E(z)E(w)&=\theta_2(z/w)\phantom{^{-1}}E(w)E(z)\\
F(z)F(w)&=\theta_2(z/w)^{-1}F(w)F(z)
\end{align*}
\item[(QL4)] Let $\delta(\zeta)=\sum_{n\in \Z}\zeta^n$ be the formal delta
function, then
\[(q-q^{-1})[E(z),F(w)]=\delta(z/w)\lp \frac{\Theta_1^+(z)}{\Theta_2^+(z)} - 
\frac{\Theta_1^-(z)}{\Theta_2^-(z)}\rp\]\\
\end{itemize}

\subsection{} Set
\begin{equation}\label{eq: psi-series}
\psi^{\pm}(z)
=\Th_1^{\pm}(z)/\Th_2^{\pm}(z)\\
=K^{\pm 1} \exp\lp\pm (q-q^{-1})\sum_{r\geq 1} 
H_{\pm r} z^{\mp r}\rp
\end{equation}
where $K=q^{H_0}$ and $H_r=D_{1,r} - D_{2,r}$, $r\in \Z$. Then,
the relation (QL4) reads
\begin{equation}\label{eq: commutation-E-F}
[E_k, F_l]=\frac{\psi^+_{k+l} - \psi^-_{k+l}}{q-q^{-1}}
\end{equation}
where $\psi^+_{-p}=\psi^-_p=0$ for every $p>0$.

\subsection{}

\begin{lem}\label{lem: commutation-nodes}
The relation (QL2) can be equivalently written as follows. For every
$r,k\in\Z$, $r\neq 0$, we have
\begin{alignat*}{3}
[D_{j,0},E_k]	&=(-1)^{j-1} E_k
&\qquad
[D_{j,0}, F_k]&=(-1)^j F_k\\
[D_{1,r}, E_k]	&= \phantom{-}q^{-r}\frac{[r]}{r}E_{k+r}
&\qquad
[D_{2,r}, E_k]	&= -q^r\frac{[r]}{r}E_{k+r}\\
[D_{1,r}, F_k] &= -q^{-r}\frac{[r]}{r}F_{k+r}
&\qquad
[D_{2,r}, F_k] &= \phantom{-}q^r\frac{[r]}{r}F_{k+r}
\end{alignat*}
And hence we have the following commutation relations
\begin{alignat*}{3}
[H_0, E_k]&=2E_k
&\qquad
[H_0,F_k]&=-2F_k\\
[H_r, E_k]&=\frac{[2r]}{r}E_{r+k}
&\qquad
[H_r, F_k]&=-\frac{[2r]}{r}F_{r+k}
\end{alignat*}
\end{lem}

\subsection{Quantum determinant}\label{sec: q-det}

It follows from the relations (QL1)--(QL2) that the coefficients of
the series
\[\qdet^\pm(z)=\Theta_1^{\pm}(q^{-1}z)\Theta_2^{\pm}(qz)\]
belong to  the center of $\hloopgl{2}$. The following result in
well--known (see \cite[Thm. 1.8.2]{molev-yangian} for the
analogous assertion for the Yangian $Y(\gl_{2})$).

\begin{prop}\label{prop: gl-to-sl}
The coefficients of $\qdet^{\pm}(z)$ generate the center $\cZ$ of
$\hloopgl{2}$. Moreover, if $\hloopsl{2}\subset\hloopgl{2}$ is the
subalgebra generated by $\{E_r, F_r,H_r\}_{r\in\Z}$, then
\[\hloopgl{2} \cong \cZ \otimes \hloopsl{2}\]
\end{prop}

\subsection{Hopf algebra structure}\label{sec: comultiplication-qla}

The algebra $\hloopgl{2}$ is a Hopf algebra with comultiplication 
determined by
\begin{equation}\label{eq: comultiplication-qla}
\begin{array}{rcl}
\Delta\lp\qdet^{\pm}(z)\rp &=& \qdet^{\pm}(z)\otimes \qdet^{\pm}(z)\\
\Delta(D_{j,0}) &=& D_{j,0}\otimes 1 + 1\otimes D_{j,0} \\
\Delta(E_0) &=& E_0\otimes K + 1\otimes E_0 \\
\Delta(F_0) &=& F_0\otimes 1 + K^{-1}\otimes F_0 \\
\Delta(E_{-1}) &=& E_{-1}\otimes K^{-1} + 1\otimes E_{-1}\\
\Delta(F_1) &=& F_1\otimes 1 + K\otimes F_1
\end{array}
\end{equation}

\subsection{Evaluation homomorphism}\label{sec: ev-qla}
For any invertible element $\zeta\in\C[[\hbar]]$, there is a surjective
algebra homomorphism $\ev_{\zeta}:\hloopgl{2}\to\Uhgl{2}$ given
on $\hloopsl{2}$ by \cite[\S 4.1]{chari-pressley-qaffine}
\begin{gather*}
H_0 \mapsto D_1-D_2
\qquad\qquad E_0 \mapsto E
\qquad\qquad F_0 \mapsto F\\[1.1 ex]
E_{-1} \mapsto q\zeta^{-1}K^{-1}E
\qquad\qquad F_1 \mapsto q^{-1}\zeta FK
\end{gather*}
and on the center $\cZ$ by 
\[\qdet^{\pm}(z)\mapsto
q^I\,\frac{z-q^{-I}\zeta}{z-q^I\zeta}\]
where $I=D_1 + D_2\in\Uhgl{2}$, and the \rhs is expanded in powers
of $z^{\mp 1}$.

If $\V$ is a $\Uhgl{2}$--module, we denote the $\hloopgl{2}$--module
$\ev_{\zeta}^*\lp\V\rp$ by $\V(\zeta)$.

\subsection{Kac--Moody presentation of $\hloopsl{2}$}\label{sec: KM}

The quantum loop algebra $\hloopsl{2}$ is usually presented on
\KM generators $\HH,\E_i,\F_i$, $i=0,1$, satisfying the relations
\begin{itemize}
\item[(KM1)] 
$\ds{[\HH, \E_i]=(-1)^{i+1}2\E_i \aand [\HH, \F_i]=(-1)^i 2\F_i}$\\
\item[(KM2)] For any $i,j\in \{0,1\}$
\[[\E_i, \F_j]=\delta_{ij} (-1)^{i+1}\,\frac{q^{\HH} - q^{-\HH}}{q-q^{-1}}\]
\item[(KM3)] For any $i\neq j \in \{0,1\}$
\begin{gather*}
\E_i^3\E_j - [3] \E_i^2\E_j\E_i + [3] \E_i\E_j\E_i^2-\E_j\E_i^3=0\\
\F_i^3\F_j - [3] \F_i^2\F_j\F_i + [3] \F_i\F_j\F_i^2-\F_j\F_i^3=0
\end{gather*}
\end{itemize}

The fact that the relations of Section \ref{sec: defn-qla} give an
equivalent presentation to the ones above is stated in \cite
{drinfeld-yangian-qaffine} (see \cite{beck-braid} for a proof). The
relation between the generators is given by\footnote{we follow
here the conventions of \cite{beck-braid}.}
\begin{equation}\label{eq: loop-to-KM}
\begin{aligned}
	&			&\HH&=D_{1,0} - D_{2,0}	&		&\\
\E_1	&=E_0		&	&				&\F_1	&=F_0 \\
\E_0	&=K^{-1}F_1	&	&				&\F_0	&=E_{-1}K
\end{aligned}
\end{equation}

\subsection{Diagram automorphism}\label{ss:diagram}

Let $\omega$ be the diagram automorphism of $\hloopsl{2}$ given
by $\HH \leftrightarrow-\HH$, $\E_0\leftrightarrow\E_1$  and $\F_0
\leftrightarrow\F_1$. In terms of loop generators, \eqref{eq: loop-to-KM}
shows that $\omega$ is given by
\begin{equation}\label{defn: diagram-automorphism}
\begin{aligned}
\omega(E_0)	&=K^{-1}F_1 	&\qquad\omega(F_0)&=E_{-1}K \\
\omega(E_{-1})&=F_0K		&\qquad\omega(F_1)&=K^{-1}E_0
\end{aligned}
\end{equation}

\section{Quantum Weyl groups}\label{sec: qWeyl-gl}

In this section, we extend the action of the affine braid group
$B_{SL_2}$ on the quantum loop algebra $\hloopsl{2}$ to one
of $B_{GL_2}$ on $\hloopgl{2}$.
We show that this action is given by conjugating by elements
$\qS,\qL_1,\qL_2$, where $\qS$ is the quantum Weyl group
element of $\Uhsl{2}$ and $\qL_1,\qL_2$ lie in a completion
of the commutative subalgebra of $\hloopgl{2}$ generated by
the elements $\{D_{j,k}\}$. The element $\qL=\qL_1\qL_2^{-1}$
is equal to the quantum Weyl group element of $\hloopsl{2}$
corresponding to the generator of the coroot lattice, for which
we obtain an explicit expression in terms of the commuting
generators $H_k$.

\subsection{Braid group action on $\hloopsl{2}$}
\label{sec: qWeyl-sl}

Following \cite{lusztig-book}, consider the automorphisms $T_0,
T_1$ of $\hloopsl{2}$ given in the \KM presentation by
\begin{gather*}
T_0(\HH)=-\HH=T_1(\HH)\\
T_i(\E_i)=-\F_i\K_i
\qquad\qquad
T_i(\F_i)=-\K_i^{-1}\E_i
\end{gather*}
where $\K_0=q^{-\HH}$, $\K_1=q^{\HH}$ and, for $i\neq j\in\{0,1\}$,
\[T_i(\E_j) = \E_i^{(2)}\E_j - q^{-1}\E_i\E_j\E_i + q^{-2}\E_j\E_i^{(2)}
\qquad
T_i(\F_j) = \F_j\F_i^{(2)} - q\F_i\F_j\F_i + q^2\F_i^{(2)}\F_j\]
where $X^{(n)}=X^n/[n]!$.

It is clear that the diagram automorphism $\omega$ defined in
\ref{ss:diagram} satisfies
\begin{equation}\label{cor: conjugation}
\omega \circ T_0 \circ \omega=T_1
\end{equation}

\subsection{} 

We shall need for later use the following
\begin{lem}\label{le:T0-loop}
The action of $T_0$ on the loop generators is given by
\[\begin{aligned}
T_0\lp F_1\rp &= -K^{-1}E_{-1} &\qquad  	T_0\lp E_{-1}\rp&= -F_1K\\
T_0\lp E_0\rp &= -K^{-1}F_2	&\qquad 	T_0\lp F_0\rp 	&= -E_{-2}K
\end{aligned}\]
\end{lem}
\begin{pf}
The first set of equations is a direct consequence of \ref{sec: qWeyl-sl}
and \eqref{eq: loop-to-KM}. We only check the first of the remaining two
equations since the second one is verified in a similar way. We have
\[T_0(E_0)=
\frac{1}{[2]}(
K^{-1}F_1K^{-1}F_1E_0
-(1+q^{-2}) K^{-1}F_1E_0K^{-1}F_1+
q^{-2} E_0K^{-1}F_1K^{-1}F_1)\]
We now rewrite each term individually.
\[K^{-1}F_1K^{-1}F_1E_0=
q^{-2}K^{-2}F_1^2E_0=q^{-2}K^{-2}(F_1E_0F_1-F_1KH_1)\]
where we used $[E_0,F_1]=\Psi^+_1/(q-q^{-1})=KH_1$. Next,
\[(1+q^{-2}) K^{-1}F_1E_0K^{-1}F_1= (1+q^{-2})K^{-2}F_1E_0F_1\]
Finally,
\[q^{-2} E_0K^{-1}F_1K^{-1}F_1=
K^{-2}E_0F_1F_1=
K^{-2}(F_1E_0F_1+KH_1F_1)\]
Combining these computations we get
\[T_0(E_0)=\frac{1}{[2]}K^{-1}[H_1,F_1]=-K^{-1}F_2\]
\end{pf}

\subsection{} 

Let $\omega\in\Aut(\hloopsl{2})$ be the diagram automorphism
defined in \S \ref{ss:diagram}.

\begin{lem}\label{le:T0w}
The action of $T_0\omega$ on $\hloopsl{2}$ is given in the
loop generators by
\[\psi^{\pm}(z)\mapsto\psi^{\pm}(z),\qquad
E(z)\mapsto -z^{-1}E(z)\aand
F(z)\mapsto -zF(z)\]
\end{lem}
\begin{pf}
It suffices to verify the assertion on the generators $K,E_0,E_{-1},F_0,
F_1$. It is clear that $T_0\omega$ fixes $K$. Moreover, using \S
\ref{ss:diagram} and Lemma \ref{le:T0-loop}, we find
\begin{alignat*}{3}
T_0(\omega(E_0))&=T_0(K^{-1}F_1)=-E_{-1}
&\qquad
T_0(\omega(E_{-1}))&=T_0(F_0K)=-E_{-2}\\
T_0(\omega(F_0))&=T_0(E_{-1}K)=-F_1
&\qquad
T_0(\omega(F_1))&=T_0(K^{-1}E_0)=-F_2
\end{alignat*}
\end{pf}

\subsection{The lattice element $L$}\label{sec: lattice-sl}

Let $L=T_0T_1\in\Aut(\hloopsl{2})$. Note that, by \eqref{cor: conjugation}
\[L=T_0T_1=T_0\omega T_0\omega=(T_0\omega)^2\]
By Lemma \ref{le:T0w}, the action of $L$ on $\hloopsl{2}$ is therefore
given by
\begin{equation}\label{eq: L-loop}
\psi^{\pm}(z)\mapsto\psi^{\pm}(z),\qquad
E(z)\mapsto z^{-2}E(z)\aand
F(z)\mapsto z^2 F(z)
\end{equation}

\subsection{The automorphisms $\mathbf{L_1,L_2}$}

Consider the assignments
\begin{equation}\label{eq:Li}
\begin{aligned}
L_1:\,
&\Theta_j^\pm(z)\to\Theta_j^\pm(z)
&\qquad
E(z)&\to z^{-1}E(z)
&\qquad
F(z)&\to z\phantom{^{-1}} F(z)\\
L_2:\,&\Theta_j^\pm(z)\to\Theta_j^\pm(z)
&\qquad
E(z)&\to z\phantom{^{-1}}E(z)
&\qquad
F(z)&\to z^{-1} F(z)
\end{aligned}
\end{equation}

\begin{prop}\label{pr:extend}\hfill
\begin{enumerate}
\item $L_1$ and $L_2$ extend uniquely to algebra automorphisms of $\hloopgl{2}$
satisfying $L_1L_2=L_2L_1=\id$.
\item The automorphism $L=T_0T_1$ is equal to $L_1L_2^{-1}$.
\item $L_1$ and $L_2$ satisy
\[T_1 L_2 T_1=L_1\]
and therefore give rise to an action of the affine braid group $B_{GL_2}$ on $\hloopgl
{2}$ extending that of $B_{SL_2}$ on $\hloopsl{2}$.
\end{enumerate}
\end{prop}
\begin{pf}
(1) It is clear from \eqref{eq:Li} that $L_1$ and $L_2$ preserve the defining
relations (QL1)--(QL4) of $\hloopgl{2}$ and that $L_1L_2=L_2L_1=\id$.
(2) Follows by comparing \eqref{eq:Li} and \eqref{eq: L-loop}.
(3) It readily follows from Lemma \ref{le:T0-loop} and \eqref{eq:Li} that $T_0L_2
T_0=L_1$.  Since $T_0T_1=L=L_1L_2^{-1}$, we have
\[L_1=T_0L_2T_0=
L_1L_2^{-1}T_1^{-1}L_2L_1L_2^{-1}T_1^{-1}=
L_1(L_2^{-1}T_1^{-1}L_1T_1^{-1})\]
Simplifying $L_1$ yields the claimed identity.
\end{pf}

\begin{rem}\label{rk:T0w Li}
Comparing Lemma \ref{le:T0w} and \eqref{eq:Li} shows that the restriction
of $L_1$ to $\hloopsl{2}$ satisfies
\[L_1=\Ad((-1)^{\HH/2})\,T_0\,\omega\]
\end{rem}

\subsection{The Quantum Weyl group of $\hloopsl{2}$ \cite{lusztig-book}}

The automorphisms $T_0,T_1$ are almost inner. Specifically, if $\wh{\hloopsl{2}}$
is the completion of $\hloopsl{2}$ \wrt its \fd representations, there are elements
$\qS_0,\qS_1$ in $\wh{\hloopsl{2}}$ such that conjugation by $\qS_i$ preserves
$\hloopsl{2}$ and is given by the automorphism $T_i$. The elements $\qS_i$ are
given by
\begin{equation}\label{eq: triple-exponential}
\begin{aligned}
\qS_0 &= \exp_{q^{-1}}\lp q^{-1}\E_0\phantom{^{-}}q^{\HH}\rp \exp_{q^{-1}}\lp 
-\F_0\rp \exp_{q^{-1}}\lp q\E_0q^{-\HH}\rp q^{\frac{\HH(\HH-1)}{2}}\\
\qS_1 &= \exp_{q^{-1}}\lp q^{-1}\E_1 q^{-\HH}\rp \exp_{q^{-1}}\lp 
-\F_1\rp \exp_{q^{-1}}\lp q\E_1\phantom{^{-}}q^{\HH}\rp q^{\frac{\HH(\HH+1)}{2}}
\end{aligned}
\end{equation}
where the $q$--exponential is defined by
\[\exp_q(x)=\sum_{n\geq 0} q^{\frac{n(n-1)}{2}} \frac{x^n}{[n]!}\]

\subsection{Completions}\label{sec: completions}

We will similarly show that the automorphisms $L_1,L_2$ are almost
inner. We begin by defining appropriate completions of $\hloopsl{2}$
and $\hloopgl{2}$.

For $\g=\sl_2$ or $\gl_2$, the completion $\wh{\hloop}$ defined below
is a flat deformation of the completion $\wh{U(\g[z,z^{-1}])}$ of the
classical loop algebra \wrt the descending chain of
ideals $J_n=U((z-1)^n\g[z,z^{-1}])$, $n\geq 0$ (see \cite[Prop. 6.3]
{sachin-valerio-1} for $\g=\sl_{2}$). For $\g=\sl_2$, $J_n$ is the $n$th
power of $J_1$ since $\g=[\g,\g]$, and $\wh{\hloopsl{2}}$ is correspondingly
defined as the completion
\[\wh{\hloopsl{2}}=\lim_{\longleftarrow}\hloopsl{2}/\J^n\]
\wrt the kernel $\J$ of the composition
\[\hloopsl{2}
\xrightarrow{\hbar\to 0}U(\sl_2[z,z^{-1}])
\xrightarrow{z\to 1}U\sl_2\]

For $\g=\gl_2$, the powers of the ideal $J_1$ are too small\footnote
{for example, $I\otimes(z-1)^2\notin\bigcup_{n\geq 1}J_1^n$, where
$I=E_{11}+E_{22}$ is the identity matrix.} and the above construction
needs to be modified as follows. For each $r\geq 0$, $t\in \Z$ and
$X=E,F,\Th_1$ or $\Th_2$, consider the element
\[X_{r;t}=\sum_{s=0}^r (-1)^s \cbin{r}{s} X_{s+t}\]
where $\Theta_{i,l}=(\Theta_{i,l}^+-\Theta_{i,l}^{-})/(q-q^{-1})$. Note
that $X_{r;t}=x\otimes z^t (1-z)^r\mod\hbar$ where $x\in\g$ is such
that $X=x\mod\hbar$. Let $\K_r$ be the two--sided ideal of $\hloopgl{2}$
generated by the elements $\{X_{r';t}\}_{r'\geq r, t\in\Z}$, and $\hbar$
if $r=1$. Finally, let $\J_n\subset\hloopgl{2}$ be the ideal
\[\J_n=
\sum_{\substack{n_1,\ldots,n_k\geq 1\\n_1+\cdots+n_k=n}}\K_{n_1}\cdots\K_{n_k}\]
Then, $\J_n$ is a descending filtration, $\J_n\J_m\subset\J_{n+m}$,
and the completion
\[\wh{\hloopgl{2}}=\lim_{\longleftarrow}\hloopgl{2}/\J_n\]
is a flat deformation of $\wh{U(L\gl_2)}$.

\begin{rem}
Note that $\hbar\in\K_1$ implies that $\hbar\J_n\subset\J_{n+1}$
for any $n\geq 0$.
\end{rem}

\subsection{}\label{ss:center}

\begin{prop}\label{prop: center-Lsl}\hfill
\begin{enumerate}
\item The center of $\wh{\hloopsl{2}}$ is trivial.
\item The center of $\wh{\hloopgl{2}}$ is generated by the elements
$\zed_r=q^rD_{1,r}+q^{-r}D_{2,r}$, $r\in\N$.
\end{enumerate}
\end{prop}
\begin{pf}
(1) follows from the fact that $\wh{\hloopsl{2}}$ is a flat deformation
of $\wh{U(L\sl_2)}$ and that the latter algebra has trivial center.

(2) By definition of the series $\qdet^{\pm}(z)$,
\begin{equation}\label{eq:qdet}
\qdet^\pm(z)
=
\Th_1^+(q^{-1}z)\Th^+_2(qz)
=
q^{\pm(D_{1,0}+D_{2,0})}\exp\lp\pm(q-q^{-1})\sum_{r\geq 1}\zed_r z^{-r}\rp
\end{equation}
Thus, the center $Z(\hloopgl{2})$
is generated by the elements $\zed_r$, $r\in\Z$. The fact that its completion
is generated by the $\zed_r$, $r\in\N$ follows from the analogous statement
in the classical case.
\end{pf}

\subsection{}

The following is straighforward and will henceforth be used implicitly.

\begin{lem}\hfill
If $\zeta\in 1+\hbar\C[[\hbar]]$, the evaluation homomorphism
$\ev_\zeta$ extends to $\wh{\hloopgl{2}}\to\Uhgl{2}$.
\end{lem}

\subsection{The operators $\qL_1,\qL_2$}
\label{sec: defn-lattice}

Define, for any $r\geq 0$ and $i=1,2$
\begin{equation}\label{eq: log-generators-gl}
\begin{aligned}
\wt{D}_{i,r}&=
D_{i,0} + \sum_{s=1}^r (-1)^s \cbin{r}{s}\frac{s}{[s]}D_{i,s}\\
\wt{H}_r&=
H_0\phantom{_{i,}}+ \sum_{s=1}^r (-1)^s \cbin{r}{s}\frac{s}{[s]} H_s
\end{aligned}
\end{equation}

\noindent The proof of the following result is given in Appendix \ref{app:J^n}.

\begin{prop}\label{prop: log-completion}
The elements $\wt{D}_{1,r},\wt{D}_{2,r}$ lie in $\J_r$. Similarly, $\wt{H}_r\in
\J^r$ for any $r\in\N$.
\end{prop}

\noindent Thus, the following are well defined elements of $\wh{\hloopgl{2}}$
and $\wh{\hloopsl{2}}$ respectively.
\begin{gather}\label{defn: defn-lattice}
\qL_1=q^{-D_1}
\exp\lp \sum_{r\geq 1} \frac{\wt{D}_{1,r}}{r} \rp 
\qquad\qquad
\qL_2=q^{-D_1}
\exp\lp \sum_{r\geq 1} \frac{\wt{D}_{2,r}}{r}\rp \\
\qL=\exp\lp \sum_{r\geq 1} \frac{\wt{H}_r}{r} \rp=
\qL_1\qL_2^{-1}
\end{gather}

\subsection{}
\label{sec: lattice-conjugation}

\begin{prop}\label{pr:L qL}
The following holds on $\wh{\hloopgl{2}}$ for $i=1,2$
\[L_i=\Ad(\qL_i)\aand L=\Ad(\qL)\]
\end{prop}
\begin{pf}
It is clear that $\Ad(\qL_i)$ fix $\Theta_j^\pm(z)$. By Lemma
\ref{lem: commutation-nodes}
\[[\wt{D}_{1,r},E(z)]=(1-q^{-1}z)^r E(z)\aand
[\wt{D}_{2,r},E(z)]=-(1-qz)^r E(z)\]
which shows that conjugation with $\ol{\qL}_i=\exp\lp 
\sum_{r\geq 1} \frac{\wt{D}_{i,r}}{r}\rp$ is given by
\[\Ad\lp\ol{\qL}_1\rp E(z)=qz^{-1}E(z)\aand
\Ad\lp\ol{\qL}_2\rp E(z)=qzE(z)\]
The relations $\Ad(\qL_1)E(z)=z^{-1}E(z)$ and $\Ad(\qL
_2)E(z)=zE(z)$ follow from this computation and the
fact that $[D_1,E(z)]=E(z)$.
The remaining relations are proved analogously.
\end{pf}

\begin{cor}\label{cor: lattice-sl}
The product $\qS_0\qS_1$ is equal to $\qL$.
\end{cor}
\begin{pf}
Propositions \ref{pr:extend} and \ref{pr:L qL} imply that
\[\Ad(\qS_0\qS_1)=T_0T_1=L_1L_2^{-1}=\Ad(\qL)\]
The result now follows from the fact that $\wh{\hloopsl{2}}$
has trivial center by Proposition \ref{ss:center}.
\end{pf}

\subsection{}

Let $\X_+,\X_-\subset\hloopgl{2}$ be the left ideals generated
by $\{E_k\}_{k\in\Z}$ and $\{F_k\}_{k\in\Z}$ respectively.

\begin{prop}\label{prop: group-like-L}
The elements $\qL_1,\qL_2,\qL$ are grouplike modulo the
subspace
\[\calN=\X_+\otimes\X_- +\X_-\otimes\X_+\]
\end{prop}
\begin{pf}
Let $\zed_r\in Z(\hloopgl{2})$ be the elements defined by \eqref
{eq:qdet}. Since $\qdet^+(z)$ is grouplike, the elements $\zed_r$
are primitive. By \cite[Prop. 4.4 (iii)]{chari-pressley-qaffine}, the
elements $H_r=D_{1,r} - D_{2,r}$ are primitive modulo $\calN$.
The same therefore holds for $D_{j,r}$ and hence for $\wt{D}_{j,r}$,
which implies the desired assertion.
\end{pf}

\begin{rem}
An alternative proof of Proposition \ref{prop: group-like-L} for
the element $\qL$ can be obtained using the results of \cite
{kirillov-reshetikhin,  lusztig-book}. Recall that for a symmetrisable
Kac--Moody algebra $\g(\bfA)$ and a node $i$ of its Dynkin diagram,
the corresponding quantum Weyl group element $\qS_i$ satisfies
\[
\Delta(\qS_i)=\RR_{i,0}^{21}\lp \qS_i\otimes \qS_i\rp
\]
where $\RR_{i,0}$ is the truncated $R$--matrix of $\Uhsl{2}^{(i)}
\subset \Uhg(\bfA)$.
Thus the elements $\qS_i$ are group--like modulo $\mathcal{N}$
and hence the same is true for $\qL=\qS_0\qS_1$.
\end{rem}

\subsection{Action on \hw vectors}\label{sec: hw-action}

For any $\lambda\in\N$, let $\V_\lambda$ be the $(\lambda+1)
$--dimensional, indecomposable representation of $\Uhgl{2}$
with highest weight vector $\Omega_\lambda$ such that $D_1
\Omega_\lambda=\lambda\Omega_\lambda$ and $D_2\Omega
_\lambda=0$. Let $\ol{\Omega}_\lambda$ be its lowest weight
vector, and $\V_\lambda(\zeta)$ the corresponding evaluation
representation of $\hloopgl{2}$, where $\zeta\in 1+\hbar\C[[\hbar]]$.
We compute below the action of the operators $\qL_1,\qL_2$
on the highest and lowest weight vectors
\[\Omega=
\bigotimes_{i=1}^k \Omega_{\lambda_i}
\in
\V_{\lambda_1}(\zeta_1)\otimes \cdots \otimes \V_{\lambda_k}(\zeta_k)
\ni
\bigotimes_{i=1}^k \ol{\Omega}_{\lambda_i}=
\ol{\Omega}\]
of a tensor product of these evaluation representations. We shall
need the following

\begin{lem}\label{le:Theta on Omega}
The following holds on $\V_\lambda(\zeta)$
\begin{alignat*}{3}
\Th_1^\pm(z)\,\Omega_\lambda
&=
q^\lambda
\frac{z-q^{-\lambda-1}\zeta}{z-q^{\lambda-1}\zeta}\,
\Omega_\lambda
&\qquad
\Th_2^\pm(z)\,\Omega_\lambda
&=\Omega_\lambda\\
\Th_1^\pm(z)\,\ol{\Omega}_\lambda
&=\ol{\Omega}_\lambda
&\qquad
\Th_2^\pm(z)\,\ol{\Omega}_\lambda
&=
q^\lambda
\frac{z-q^{-\lambda+1}\zeta}{z-q^{\lambda+1}\zeta}\,
\ol{\Omega}_\lambda
\end{alignat*}
\end{lem}
\begin{pf}
By \cite[\S 4.2]{chari-pressley-qaffine},
\[\psi^\pm(z)\,\Omega_\lambda
=
q^\lambda
\frac{z-q^{-\lambda-1}\zeta}
{z-q^{\lambda-1}\zeta}\,\Omega_\lambda
\aand
\psi^\pm(z)\,\ol{\Omega}_\lambda
=
q^{-\lambda}
\frac{z-q^{\lambda+1}\zeta}
{z-q^{-\lambda+1}\zeta}\,\ol{\Omega}_\lambda\]
By \S \ref{sec: ev-qla}, the series $\qdet^\pm(z)$ acts on $\V_\lambda
(\zeta)$ as multiplication by
$q^\lambda\,\frac{z-q^{-\lambda}\zeta}{z-q^\lambda\zeta}$.
Using $\Psi^\pm(z)=\Th_1^\pm(z)/\Th^\pm_2(z)$ and $\qdet^\pm(z)=
\Theta_1^\pm(q^{-1}z)\Theta_2^\pm(qz)$ shows that
\[\Theta_2^\pm(qz)\Theta_2^\pm(q^{-1}z)\Omega_\lambda=
\Omega_\lambda
\aand
\Theta_1^\pm(qz)\Theta_1^\pm(q^{-1}z)\ol{\Omega}_\lambda=
\ol{\Omega}_\lambda\]
from which the stated formulae readily follow.
\end{pf}

\begin{prop}\label{prop: hw-action-L}
The following holds on $\V_{\lambda_1}(\zeta_1)\otimes \cdots
\otimes\V_{\lambda_k}(\zeta_k)$
\begin{alignat*}{3}
\qL_1\,\Omega&=\prod_{1\leq a\leq k}
\zeta_a^{-\lambda_a}\,\Omega
&\qquad
\qL_2\,\Omega&=\prod_{1\leq a\leq k}
q^{-\lambda_a}\,\Omega\\[1.1 ex]
\qL_1\,\ol{\Omega}&=\ol{\Omega}
&\qquad
\qL_2\,\ol{\Omega}&=\prod_{1\leq a\leq k}
q^{-\lambda_a}\zeta_a^{-\lambda_a}\,\ol{\Omega}
\end{alignat*}
\end{prop}
\begin{pf}
By Proposition \ref{prop: group-like-L}, it suffices to prove the
result for $k=1$. We consequently drop the subscript $a$ from
the computations below. By Lemma \ref{le:Theta on Omega},
$D_{2,r}\Omega=0$ for any $r\geq 0$, $D_{1,0}\Omega=
\lambda\Omega$ and, for $r\geq 1$, 
\[D_{1,r}\Omega=q^{-r}\frac{[\lambda r]}{r} \zeta^r \Omega\]
This implies that $\wt{D}_{2,r}\Omega=0$ and
\[\wt{D}_{1,r}\Omega=\lp\sum_{t=0}^{\lambda-1} \lp 1-\zeta
q^{2t-\lambda}\rp^r \rp\Omega\]
Thus, $\qL_2\Omega=q^{-D_1}\Omega=q^{-\lambda}\,\Omega$
and $\qL_1\Omega=q^{-D_1}(\zeta^{-\lambda}q^\lambda)\Omega=
\zeta^{-\lambda}\Omega$ as claimed. The remaining relations
follows similarly.
\end{pf}

\subsection{Braid relations}\label{sec: braid-relns}

We now show that the operators $\qS_1,\qL_1,\qL_2$ satisfy
relations very similar to those defining the affine braid group
$B_{GL_2}$.

\begin{lem}\label{lem: separation-points}
Let $\V=\V_1$ be the standard two--dimensional representation
of $\Uhgl{2}$. Then, the evaluation representations
\begin{equation}\label{eq:eval reps}
\V(\ul{\zeta})=\V(\zeta_1)\otimes\cdots\otimes \V(\zeta_k)
\end{equation}
separate the elements of the centre of $\wh{\hloopgl{2}}$.
\end{lem}
\begin{pf} It follows from \S \ref{sec: ev-qla}, and the fact that
the series $\qdet^\pm(z)$ are grouplike that their action on
$\V(\zeta)$ is given by multiplication by
\[q^{\pm k}\prod_{a=1}^k
\lp \frac{1 - q^{\mp 1} \zeta_a^{\pm 1}
z^{\mp 1}}{1-q^{\pm 1}\zeta_a^{\pm 1}z^{\mp 1}}\rp\]
Thus, the generators $\zed_r$, $r\in\N$ of $Z(\hloopgl{2})$
defined by \eqref{eq:qdet} act as multiplication by the power
sums
\[\zed_r=\frac{[r]}{r}\,\sum_{a=1}^k \zeta_a^r\]
The claim now follows from the fact that these are algebraically
independent.
\end{pf}

\begin{thm}\label{thm: first-main-theorem}
The elements $\qS_1,\qL_1,\qL_2$ satisfy the following relations
\begin{enumerate}
\item $\qL_1\qL_2=\qL_2\qL_1$.
\item $\qS_1\qL_2\qS_1=(-1)^I\qL_1$.
\end{enumerate}
where $I=D_{1,0} + D_{2,0}$.
\end{thm}
\begin{pf}
The first assertion is obvious since $\qL_1,\qL_2$ are defined in terms
of the commutating elements $D_{i,r}$. By Propositions \ref{pr:extend}
and \ref{pr:L qL}, both sides of (2) define the same automorphism of
$\hloopgl{2}$ and therefore agree up a central element $c$, $\qS_1\qL
_2\qS_1=c \qL_1$. To determine $c$ it suffices, by Lemma
\ref{lem: separation-points}, to compute it on all evaluation representations
\eqref{eq:eval reps}. Let $\Omega,\ol{\Omega}$ be the highest and lowest
weight vectors in $\V(\Omega)$. By \cite{lusztig-book}
\[\qS_1 \Omega =(-1)^k q^k \ol{\Omega}\aand
\qS_1 \ol{\Omega}=\Omega\]
Together with Proposition \ref{prop: hw-action-L}, this implies that
\[\qS_1\qL_2\qS_1\ol{\Omega}=
\qS\qL_2 \Omega=
q^{-k}\qS\Omega=
(-1)^k\ol{\Omega}=
(-1)^k\qL_1\ol{\Omega}\]
so that $c$ acts as $(-1)^I$ on $\V(\ul{\zeta})$ as claimed.
\end{pf}

\subsection{The quantum Weyl group of $\hloopgl{2}$}\label{ss:qW gl}

Set
\begin{equation}\label{eq:S S1}
\qS=\qS_1(-1)^{D_1}=(-1)^{D_2}\qS_1
\end{equation}
By Theorem \ref{thm: first-main-theorem},
the elements $\qS,\qL_1,\qL_2$ satisfy the defining relations of $B_
{GL_2}$, namely
\[\qL_1\qL_2=\qL_2\qL_1\aand \qS\qL_2\qS=\qL_1\]
We shall refer to $\qS,\qL_1,\qL_2$ as the {\it quantum Weyl group
elements} of $\hloopgl{2}$. Note that the fact that the element $\qS$
differs from $\qS_1$ by the sign $(-1)^{D_1}$ is in agreement with
the fact that their classical limits are, respectively
\[\begin{pmatrix}0&1\\1&0\end{pmatrix}\aand
\begin{pmatrix}\phantom{-}0&1\\-1&0\end{pmatrix}\]
which are the generators of the (Tits extensions of the) Weyl groups
of $GL_2$ and $SL_2$.

\section{The dual pair $\lp \Uhgl{k}, \Uhgl{n}\rp$}\label{sec: q-howe}

In this section we review a deformation of the matrix space $\cfun{k}
{n}$ as a joint representation space for $\Uhgl{k}$ and $\Uhgl{n}$.
The main reference for this section is \cite[\S 5]{valerio4}.

\subsection{Quantum matrix $(k\times n)$ space}\label{sec: q-space}

By definition, $\qfun{k}{n}$ is the algebra over $\C[[\hbar]]$ topologically 
generated by elements $\{X_{ai}\_{1\leq a\leq k, 1\leq i\leq n}$ subject
to the relations
\[
X_{ai}X_{bj}=\left\{ \begin{array}{ll}
X_{bj}X_{ai} & \text{if $a<b$ and $i>j$ or $a>b$ and $i<j$}\\
q^{-1}X_{bj}X_{ai} & \text{if $a=b$ and $i<j$ or $a<b$ and $i=j$}\\
X_{bj}X_{ai} - (q-q^{-1}) X_{bi}X{aj} & \text{if $a>b$ and $i>j$}
\end{array}\right.
\]
For each $m=(m_{ai})_{1\leq a\leq k; 1\leq i\leq n}$ define
\begin{align*}
X^m &= \lp X_{11}^{m_{11}}\cdots X_{k1}^{m_{k1}}\rp \cdots \lp 
X_{1n}^{m_{1n}}\cdots X_{kn}^{m_{kn}}\rp \\
&= \lp X_{11}^{m_{11}}\cdots X_{1n}^{m_{1n}}\rp \cdots \lp X_{k1}
^{m_{k1}}\cdots X_{kn}^{m_{kn}}\rp
\end{align*}
Then, the set $\{X^m\}_{m\in \M_{k\times n}(\N)}$ is a basis for $\qfun{k}
{n}$ over $\C[[\hbar]]$.

\subsection{The joint action of $\lp \Uhgl{k}, \Uhgl{n}\rp$}
\label{sec: joint-action}

Define the following operators on $\qfun{k}{n}$, for each $b\in\{1,
\ldots, k\}$ and $a\in \{1,\ldots, k-1\}$:
\begin{align*}
D_b^{(k)}X^m &= \sum_i m_{bi} X^m \\
E_a^{(k)}X^m &= \sum_{i=1}^n [m_{a+1,i}] \prod_{j=i+1}^n 
q^{(m_{aj} - m_{a+1,j})} X^{m+\ve_{ai} - \ve{a+1,i}} \\
F_a^{(k)}X^m &= \sum_{i=1}^n [m_{ai}] \prod_{j=1}^{i-1} 
q^{-(m_{aj}-m_{a+1,j})} X^{m-\ve_{aj} + \ve_{a+1,j}}
\end{align*}
Similarly define the operators for each $j\in \{1,\ldots, n\}$ and 
$i\in \{1,\ldots, n-1\}$:
\begin{align*}
D_j^{(n)}X^m &= \sum_{a=1}^k m_{aj} X^m\\
E_i^{(n)}X^m &= \sum_{a=1}^k [m_{a,i+1}] \prod_{b=a+1}^k 
q^{m_{b,i} - m_{b,i+1}} X^{m+\ve_{ai}-\ve_{a,i+1}}\\
F_i^{(n)}X^m &= \sum_{a=1}^k [m_{ai}] \prod_{b=1}^{a-1} 
q^{-(m_{b,i} - m_{b,i+1})} X^{m-\ve_{ai} + \ve_{a,i+1}}
\end{align*}

\noindent The following result is proved in \cite[Thm. 5.4]{valerio4}
and builds upon the approach to quantum matrix space described
in \cite{baumann}.

\begin{thm}
The operators above define a structure of an algebra module
on $\qfun{k}{n}$ over $\Uhgl{k}\otimes \Uhgl{n}$. Moreover as
a $\Uhgl{k}$ (resp. $\Uhgl{n}$) module we have
\[
\qfun{k}{n} \cong \qfun{k}{1}^{\otimes n} \ \lp\text{resp. } \qfun{1}
{n}^{\otimes k}\rp
\]
\end{thm}

\section{Affine braid group actions on quantum matrix space}\label{sec: equivalence}

\subsection{}\label{sec: statement-second}

In this section, we compare two actions of the affine braid group
$B_{GL_2}$ on the quantum matrix space $\qfun{k}{2}$ described
in Section \ref{sec: q-howe}. The first is described in \ref{ss:EG n=2}
and arises by regarding $\qfun{k}{2}$ as the $\Uhgl{k}$--module
$\qfun{k}{1}^{\otimes 2}$. It is given in the generators $b,\LL_1,
\LL_2$ of Section \ref{sec: pi-2} by
\[b\mapsto(1\,2)\RR\qquad 			
\LL_1\mapsto(q^{2\fs})^{(1)}\RR\qquad 	
\LL_2\mapsto\RR^{-1}(q^{2\fs})^{(2)}\]	
and depends upon the choice of a diagonal matrix
\[s=\sum_{a=1}^ks_a E_{aa}\in\gl_k\]

The second action is obtained by regarding $\qfun{k}{2}$ as the
tensor product of evaluation representations of $\hloopgl{2}$
\[\qfun{k}{2}\cong
\qfun{1}{2}(\zeta_1)\otimes\cdots\otimes\qfun{1}{2}(\zeta_k)\]
corresponding to a choice of evaluation points $\ul{\zeta}
=(\zeta_1,\ldots,\zeta_k)\in(1+\hbar\C[[\hbar]])^k$. It is
given in terms of the quantum Weyl group elements of
$\hloopgl{2}$ defined in \ref{ss:qW gl} by
\[b\mapsto\qS\qquad\qquad\LL_1\mapsto\qL_1\qquad\qquad\LL_2\mapsto\qL_2\]

It was shown in \cite{valerio4} that the restrictions of these actions
to the braid group $B\subset B_{GL_2}$ generated by $b$ essentially
coincide. Specifically, one has $(12)\RR_{\Lsl_k}=\qS_1 q^{-(D_1+
D_1D_2/k)}(-1)^{D_1}$ \cite[Thm. 6.5]{valerio4} which, by Remark
\ref{rem: R-matrix-comparison} and \eqref{eq:S S1}, implies that
\begin{equation}\label{eq:R=S}
(1\,2)\RR=\qS\,q^{-D_1}
\end{equation}
The result below shows that similar relations hold between the
operators giving the actions of the generators $\LL_1,\LL_2$.

\begin{thm}\label{thm: second-main-theorem}
Assume that the evaluation points $\zeta_1,\ldots,\zeta_k$ are given by
\begin{equation}\label{eq:zeta q}
\zeta_a=q^{-2s_a}
\end{equation}
Then, the following holds on $\qfun{k}{2}$
\begin{align}
(q^{2\fs})^{(1)}\RR&=\qL_1					\label{eq:L1=L1}\\
\RR^{-1}(q^{2\fs})^{(2)}&=\qL_2\,q^I				\label{eq:L2=L2}
\end{align}
\end{thm}
\begin{pf}
It is easy to see, using $\LL_2=b^{-1}\LL_1b^{-1}$ and $\qL_2=\qS
^{-1}\qL_1\qS^{-1}$ that \eqref{eq:R=S} and \eqref{eq:L1=L1} imply
\eqref{eq:L2=L2}. The proof of \eqref{eq:L1=L1} occupies the rest
of this section. We first show in Proposition \ref{cor: proof-b-1} that
both sides of \eqref{eq:L1=L1} have the same commutation relation
with elements in $\hloopsl{2}$. We then check in Lemma \ref{lem: coincidence-3}
that they coincide on the tensor product of highest weight vectors in
\[\qfun{1}{2}[\lambda_1]\otimes\cdots\otimes\qfun{1}{2}[\lambda_k]
\subset
\qfun{k}{2}\]
where the notation $[\lambda_i]$ refers to the homegeneity degree
in the variables $X_{i1},X_{i2}$. If the evaluation points are generic,
the statement follows because the action of $\hloopsl{2}$ on the 
above tensor product is irreducible. The general case follows by
continuity.
\end{pf}

\subsection{}\label{sec: proof-begin}

Let $\tau=(1\,2)$ be the flip acting on $\qfun{k}{2}\cong\qfun{k}{1}
^{\otimes 2}$. In terms of the monomial basis $\{X^m\}$, the action
of $\tau$ is given by
\[(X_{11}^{m_{11}}X_{12}^{m_{12}})\cdots (X_{k1}^{m_{k1}}
X_{k2}^{m_{k2}}) \mapsto (X_{11}^{m_{12}}X_{12}^{m_{11}})
\cdots (X_{k1}^{m_{k2}}X_{k2}^{m_{k1}})\]

\begin{lem}\label{lem: coincidence-1}
The following holds on $\qfun{k}{2}$
\begin{enumerate}
\item $\ds{(q^{2\fs})^{(1)}=q^{2s_1D_1}\otimes\cdots\otimes q^{2s_kD_1}}$.\\
\item For any $x\in\Uhgl{2}^{\otimes k}$
\[\Ad\left(\tau\right)\,x=\ast^{\otimes k}(x)\]
where $\ast\in\Aut(\Uhgl{2})$ is the involution given by
\[D_1\leftrightarrow D_2\aand E\leftrightarrow F\]
\end{enumerate}
\end{lem}
\begin{pf}
(1) and (2) follow from the formulae giving the action of $\Uhgl{2}$ in
\ref{sec: joint-action}.
\end{pf}

\subsection{}  

\newcommand {\lexp}{q^{\HH/2}}

\begin{lem}\label{lem: coincidence-2}
Assume that the evaluation points for $\hloopsl{2}$ are given by 
\eqref{eq:zeta q}.
Then, the following holds on $\qfun{k}{2}$ for any $X\in\hloopsl{2}$
\[\Ad\left((q^{2\fs})^{(1)}\,\tau\right)\,X=\Ad(\lexp)\omega(X)\]
where $\omega\in\Aut(\hloopsl{2})$ is the diagram automorphism
defined in \ref{ss:diagram}.
\end{lem}
\begin{pf}
The stated identity clearly holds for $X=\HH$. It therefore suffices
to check it on the remaining generators $E_0,F_0,E_{-1},F_{-1}$
of $\hloopsl{2}$. Moreover, since
\[\left((q^{2\fs})^{(1)}\,\tau\right)^2=q^{2\fs}\otimes q^{2\fs}\]
commutes with the action of $\hloopsl{2}$, $\Ad((q^{2\fs})^{(1)}\,\tau)$
acts as an involution on the image of $\hloopsl{2}$. Since so does
$\Ad(\lexp)\omega$, it suffices to check the identity on only one half of these
generators which, in view of the formulae \eqref{defn: diagram-automorphism}
can be taken to be $E_0,F_0$.

By \eqref{eq: comultiplication-qla}, $E_0$ acts on $\qfun{k}{2}$ by
\[\ev_{\ul{\zeta}}\Delta^{(k)}(E_0)=
\sum_{a=1}^k 1^{\otimes (a-1)}\otimes E \otimes K^{\otimes (k-a)}\]
so that, by Lemma \ref{lem: coincidence-1}
\[\Ad\left((q^{2\fs})^{(1)}\,\tau\right)\ev_{\ul{\zeta}}\Delta^{(k)}(E_0)=
\sum_{a=1}^k 1^{\otimes (a-1)}\otimes q^{-2s_a}F
\otimes (K^{-1})^{\otimes (k-a)}\]
\noindent On the other hand, the $k$--fold coproduct of $\Ad(\lexp)
\omega(E_0)=q^{-1}K^{-1}F_1$ is equal to
\[(K^{-1})^{\otimes k}\,
\sum_{a=1}^k K^{\otimes (a-1)}\otimes q^{-1}F_1\otimes 1^{\otimes (k-a)}
=
\sum_{a=1}^k
1^{\otimes (a-1)}\otimes q^{-1}K^{-1}F_1 \otimes (K^{-1})^{\otimes (k-a)}\]
so that its image under $\ev_{\ul{\zeta}}$ is equal to
\[\sum_{a=1}^k
1^{\otimes (a-1)}\otimes q^{-2}\zeta_aK^{-1}FK \otimes (K^{-1})^{\otimes (k-a)}=
\sum_{a=1}^k
1^{\otimes (a-1)}\otimes \zeta_a F\otimes (K^{-1})^{\otimes (k-a)}\]
The computation for $F_0$ is identical.
\end{pf}

\subsection{}

\begin{prop}\label{cor: proof-b-1}
Assume that the evaluation points are given by \eqref{eq:zeta q}.
Then, the following holds on $\qfun{k}{2}$ for any $X\in\hloopsl{2}$
\[\Ad\left((q^{2\fs})^{(1)}\RR\right)X=\Ad(\qL_1)X\]
\end{prop}
\begin{pf}
By \eqref{eq:R=S} and Lemma \ref{lem: coincidence-2}, the \lhs
is equal to
\[\begin{split}
\Ad\left((q^{2\fs})^{(1)}\tau\,\qS_1q^{-D_1}(-1)^{D_1}\right)\,X
&=
\Ad(q^{-D_1}(-1)^{D_1}\lexp)\omega\,T_1(X)\\
&=
\Ad(q^{-I/2}(-1)^{D_1})T_0\,\omega(X)\\
&=
L_1(X)\\
&=
\Ad(\qL_1)X
\end{split}\]
where the second equality uses \eqref{cor: conjugation}, the third one
Remark \ref{rk:T0w Li}, and the last one Proposition \ref{pr:L qL}.
\end{pf}

\subsection{}

Let $\lambda=(\lambda_1,\ldots,\lambda_k)\in\N^k$ and set
\[\Omega=
X_{11}^{\lambda_1}X_{21}^{\lambda_2}\cdots X_{k1}^{\lambda_k}\in
\qfun{k}{2}\]

\begin{lem}\label{lem: coincidence-3}
The following holds
\[(q^{2\fs})^{(1)}\RR\,\Omega=\prod_{a=1}^k q^{2s_a\lambda_a}\,\Omega
\aand
\qL_1\,\Omega=\prod_{a=1}^k\zeta_a^{-\lambda_a}\,\Omega\]
\end{lem}
\begin{pf}
Under the identification $\qfun{k}{2}\cong\qfun{k}{1}^{\otimes 2}$ of
$\Uhgl{k}$--modules, $\Omega$ is the tensor product of a vector in the
$q$--deformation of $S^{\sum_a\lambda_a}\C^k$ and a vector in the
the trivial representation of $\Uhgl{k}$. Thus $\RR\Omega=\Omega$,
which implies the first stated formula. Under the identification $\qfun{k}
{2}\cong\qfun{1}{2}^{\otimes k}$ of $\Uhgl{2}$--modules, $\Omega$ is
the tensor product of highest weight vectors in $\V_{\lambda_1}\otimes
\cdots\otimes\V_{\lambda_k}$, where the notation is as in \ref{sec: hw-action}.
The result then follows from Proposition \ref{prop: hw-action-L}.
\end{pf}

\section{Monodromy theorems}\label{se:monodromy}

For any $\lambda\in\N$, denote by $V_\lambda=S^\lambda\C^2$
the $\lambda$th symmetric power of the defining representation of
$\gl_2$, and by $\V_\lambda$ its quantum deformation, that is the
\fd $\Uhgl{2}$--module such that $\V_\lambda/\hbar\V_\lambda
\cong V_\lambda$ and $I$ acts as multiplication by $\lambda$ on
$\V_\lambda$. 

Fix now $\ul{\lambda}=(\lambda_1,\ldots,\lambda_k)\in\N^k$, $\ul
{s}=(s_1,\ldots,s_k)\in\C^k$, and denote by
\[V_{\ul{\lambda}}(\ul{s})=
V_{\lambda_1}(s_1)\otimes\cdots\otimes V_{\lambda_k}(s_k)\]
the tensor product of the evaluation modules of the Yangian $\Ygl
{2}$ corresponding to the points $r_a=\fp(s_a+\frac{I+1}{2})$.
By Lemma \ref{le:shift}, the restriction of $V_{\ul{\lambda}}(\ul{s})$
to $\Ysl{2}\subset\Ygl{2}$ is the tensor product of the modules $V
_{\lambda_1},\ldots,V_{\lambda_k}$ evaluated at the points $\fp
s_1,\ldots,\fp s_k$. Denote by
\[\V_{\ul{\lambda}}(\ul{\zeta})=\V_{\lambda_1}(\zeta_1)\otimes\cdots\otimes\V_{\lambda_k}(\zeta_k)\]
the tensor product of evaluation modules of the quantum loop algebra
$\hloopgl{2}$ corresponding to the evaluation points $(\zeta_1,\ldots,
\zeta_k)\in(\C[[\hbar]]^\times)^k$. The following is the main result of
this paper.

\begin{thm}
Let $\g=\sl_2$ or $\gl_2$. Assume that $\hbar=4\pi\imath\fp$, and that 
$\zeta_a=\exp(-\hbar s_a)$ for any $a$. Then, the monodromy of the trigonometric Casimir
connection of $\g$ on $V_{\ul{\lambda}}(\ul{s})$ is described by the
quantum Weyl operators of the quantum loop algebra $\hloop$ on 
$\V_{\ul{\lambda}}(\ul{s})$.
\end{thm}
\begin{pf}
We first prove the result for $\g=\gl_2$. Let $\cfun{k}{2}$ be the space
of functions on the space of $k\times 2$ matrices described in Section
\ref{sec: howe}. As a $U\gl_2^{\otimes k}$--module, $V_{\ul{\lambda}}
(\ul{s})$ may be realised as the subspace of $\cfun{k}{2}$ via
\[V_{\lambda_1}\otimes\cdots\otimes V_{\lambda_k}\subset
S^\bullet\C^2\otimes\cdots\otimes S^\bullet\C^2\cong
\cfun{k}{2}\]
Combining the duality statement of Corollary \ref{cor: monodromy-comparison}
with the computation of the monodromy of the trigonometric KZ connection
for $\gl_k$ on $n=2$ points given in \ref{ss:EG n=2} and Theorem \ref
{thm: second-main-theorem}, we see that the monodromy of the trigonometric
connection of $\gl_2$ on $V_{\ul{\lambda}}(\ul{s})$ and the quantum Weyl
group operators $\qS,\qL_1,\qL_2$ giving the action of $\BGL$ on the
quantum matrix space $\qfun{k}{2}$ are related by
\begin{gather*}
\pi_{C,\ul{s}}(b)=
\pi_{\KZ,s}(b)q^{I/2}=(1\,2)\RR\,q^{I/2}=\qS\,q^{-\HH/2}\\
\pi_{C,\ul{s}}(\LL_1)=
\pi_{\KZ,s}(\LL_1)=(q^{2s})^{(1)}\RR=\qL_1\\
\pi_{C,\ul{s}}(\LL_2)=
\pi_{\KZ,s}(\LL_2)q^{-I}=\RR^{-1}(q^{2s})^{(2)}\,q^{-I}=\qL_2
\end{gather*}
where $b,\LL_1,\LL_2$ are the generators of $\BGL$ described in
\ref{sec: inc-ex-end}. The assertion of the theorem now follows
since $\qS\,q^{-\HH/2}=q^{\HH/4}\,\qS\,q^{-\HH/4}$. For $\g=\sl_2$,
the corresponding braid group is generated by $b,\LL=\LL_1\LL_2
^{-1}$. Moreover,
\[\pi_{C,\ul{s}}(\LL)=\pi_{C,\ul{s}}(\LL_1)\pi_{C,\ul{s}}(\LL_2)^{-1}=\qL_1\qL_2^{-1}=\qL\]
and
\[\pi_{C,\ul{s}}(b)=\pi_{C,\ul{s}}(b)(-1)^{E_{11}}=\qS(-1)^{D_1}\,q^{-\HH/2}=\qS\,q^{-\HH/2}\]
\end{pf}

\appendix

\section{Monodromy of the trigonometric KZ equations
(after Etingof--Geer--Schiffmann)}
\label{app-sec: etingof-geer}

This appendix follows \cite{etingof-geer} closely. It only differs from
it in the explicit description of the monodromy of the trigonometric
KZ equations, which is not quite correct as stated in \cite[Thm. 3.3]
{etingof-geer}.\footnote{In turn, \cite{etingof-geer} amends the
computation of the monodromy of the trigonometric KZ equations
given in \cite{etingof-schiffmann-ed1}.}

\subsection{Trigonometric KZ equations}

Let $A$ be a unital, associative algebra over $\C$ and $r\in A
\otimes A$ a classical $r$--matrix, that is a solution of the \CYBE
\begin{equation}\label{app-eq: cybe}
[r_{12}, r_{23}] + [r_{12}, r_{13}] + [r_{13},r_{23}]=0
\end{equation}
Set $r(u)=\ds{\frac{re^u+r_{21}}{e^u-1}}$, and let $s\in A$ be such
that 
\begin{equation}\label{app-eq: r-s}
[s\otimes 1+1\otimes s,r]=0
\end{equation}
Let $V$ be an $A$--module, $n\geq 1$, and $\IV^{\otimes n}$ the trivial
bundle over $\C^n$ with fibre $V^{\otimes n}$. The trigonometric KZ
connection is the flat, $\Sym_n$--equivariant connection on $\IV^{\otimes n}$
given by
\begin{equation}\label{app-eq: trig-KZ}
\nabla\KKZ=
d-\frac{\hbar}{2\pi\iota}\left(\sum_{i<j}r_{ij}(u_i-u_j)d(u_i-u_j)+\sum_i s^{(i)}du_i\right)
\end{equation}
As explained in Section \ref{ss:Pi_n}, its monodromy yields a representation
\begin{equation}\label{app-eq: monodromy}
\pi_{KZ} : \Pi_n \to GL\lp V^{\otimes n}[[\hbar]]\rp
\end{equation}
of the fundamental group $\Pi_n$ of the configuration space of $n$ points in
$\nC$ on $V^{\otimes n}[[\hbar]]$.

\subsection{\EK quantization} 

In order to describe the monodromy representation \eqref
{app-eq: monodromy}, one uses the machinery of \EK
quantization \cite{etingof-kazhdan-quantization-1}. Define
the following \fd subspaces of $A$
\begin{gather*}
\g_+=\{(1\otimes f)(r) : f\in A^*\}\\
\g_-=\{(g\otimes 1)(r) : g\in A^*\}
\end{gather*}
The following is a consequence of \eqref{app-eq: cybe} (see
\cite[\S 5]{etingof-kazhdan-quantization-1} for details)
\begin{prop}\label{app-prop: manin-triple}\hfill 
\begin{enumerate}
\item $\g_{\pm}$ are Lie subalgebras of $A$.
\item The following defines a non--degenerate pairing $\g_+
\otimes\g_-\to\C$
\[\la (1\otimes f)(r), (g\otimes 1)(r)\ra=(g\otimes f)(r)\]
\item The vector space $\g=\g_+\oplus\g_-$ is endowed with a
unique Lie algebra structure extending those on $\g_\pm$ and
such that
\[[x_+,x_-]=\ad^*(x_+)x_--\ad^*(x_-)x_+\]
where $x_\pm\in\g_\pm$ and $\ad^*$ denotes the coadjoint
action of $\g_\pm$ on $\g_\mp\cong\g_\pm^*$.
\item The map $\pi:\g\to A$ whose restriction to $\g_{\pm}$ is
the canonical inclusion is a Lie algebra homomorphism.
\item The canonical element $1\in\End(\g_+)\cong\g_+\otimes
\g_-$ maps to $r$ under the homomorphism $\pi:\g\to A$.
\end{enumerate}
\end{prop}

It follows from Proposition \ref{app-prop: manin-triple} that $(\g,
\g_+,\g_-)$ is a Manin triple. Using the quantization theorem for
\fd Manin triples \cite[\S 3]{etingof-kazhdan-quantization-1}, we
obtain a \qt Hopf algebra $\Uhg$ with $R$--matrix $\wtR\in\Uhg
^{\otimes 2}$, and Hopf subalgebras
$\Uhg_{\pm}\subset\Uhg$ in duality with each other, such that
$\wtR\in\Uhg_+\otimes\Uhg_-$. Moreover, there is a canonical
isomorphism of algebras $\Uhg\to U\g[[\hbar]]$ which allows us
to extend the map $\pi : \g \to A$ to a homomorphism $\Uhg\to
A[[\hbar]]$. 

Consider now the following elements\footnote{The elements $T,C$
differ slightly from those defined in \cite{etingof-geer} which are,
respectively, $T'=\Id\otimes S(\wtR)=T_{21}$ and
$C'=m_{01}(T_{01}\cdots T_{0n})$}
\begin{align}
T&= S\otimes\Id(\wtR_{21})\in\Uhg^{\otimes 2}\\
C&=m_{01}(T_{0n}\cdots T_{01})=m_{01}(\Id\otimes\Delta^{(n)}(T))\in\Uhg^{\otimes n}
\label{app-eq: eg-element}
\end{align}
where $m_{01}$ is the multiplication on the first two copies in $\Uhg
^{\otimes (n+1)}$, $\Delta^{(n)}:\Uhg\to\Uhg^{\otimes n}$ the iterated
coproduct, and the second equality in \eqref{app-eq: eg-element} follows
from the cabling identity $\Delta^{(n)}\otimes\Id(\wtR)=\wtR_{0\,n}\wtR
_{1\,n}\cdots\wtR_{n-1\,n}$ which implies that
\[\Id\otimes\Delta^{(n)}(\wtR_{21})=(0\,1\cdots n)\Delta^{(n)}\otimes\Id(\wtR)=
\wtR_{1\,0}\wtR_{2\,0}\cdots\wtR_{n\,0}\]
Let $\V$ the $\Uhg$--module obtained from $V\fml$
via the homomorphism $\pi:\Uhg\to A\fml$. Then, the following holds

\begin{thm}\label{app-thm: etingof-geer}
The monodromy representation $\pi_{KZ}$ \eqref{app-eq: monodromy}
is equivalent to the action of $\Pi_n$ on $\V^{\otimes n}$ given by
\begin{align*}
b_i 		&\mapsto (i\,i+1)\,\wtR_{i,i+1} \\
\X_1		&\mapsto (e^{\hbar s} u^{-1})^{(1)} C
\end{align*}
where $u=m(S\otimes\Id)(\wtR_{21})\in\Uhg$ is the Drinfeld element.
\end{thm}

The proof of this theorem is sketched in \S \ref{app-sec: proof-begin}
-- \S \ref{app-sec: proof-end}.

\subsection{}\label{app-sec: proof-begin} 

The first step towards the proof of Theorem \ref{app-thm: etingof-geer}
is to relate the trigonometric KZ connection on $n$ points \eqref
{app-eq: trig-KZ} to the rational KZ connection on $n+1$ points.
This is achieved by extending the Manin triple of Proposition \ref
{app-prop: manin-triple} to include a derivation.

Let $\rho_r=m(r_{21})\in A$, so that if $r=\sum_i a_i\otimes b_i$,
then $\rho_r=\sum_i b_ia_i$, and note that $[s,\rho_r]=0$ by \eqref
{app-eq: r-s}, Set $t=s+\rho_r$. It follows from \eqref{app-eq: cybe}
that $[t\otimes1+1\otimes t,r]=0$, which implies that $\ad(t)\g_\pm
\subset\g_\pm$, and that $\ad(t)$ is a derivation of $\g_\pm$
preserving the pairing between $\g_+$ and $\g_-$. Let $\g'=
(\g\rtimes\C t)\oplus\C t^*$ be the extension of $\g\rtimes\C t$
by a central element $t^*$ determined by requiring that the
commutator with $t$ is the derivation $\ad(t)$ on $\g$ and
that, for $x,y\in\g$
\[[x,y]_{\g'}=[x,y]_{\g}+([t,x],y)\,t^*\]
Note that $\g'$ is split over $\g_\pm\rtimes\C t$. The inner
product on $\g$ extends to a non--degenerate, invariant
bilinear form $(-,-)$ on $\g'$ given by 
\[(t,\g)=(t^*,\g)=0\aand (t,t^*)=1\]
Thus, $(\g',\g_+'=\g_+\rtimes\C t,\g'_-=\g_-\oplus\C t^*)$ is a Manin
triple. The corresponding Lie cobracket $\delta_{\g'}:\g'\to\g'\wedge\g'$
is given by $\delta_{\g'}(t)=\delta_{\g'}(t^*)=0$, $\delta_{\g'}(x)=
\delta_{\g}(x)$ if $x\in\g_+$ and
\[\delta_{\g'}(x)=\delta_\g(x)+[t,x]\wedge t^*\]
if $x\in\g_-$. In particular, $\g_+$ is a Lie subbialgebra of $\g_+$,
but $\g_-$ is only a Lie subalgebra of $\g_-'$.

Extend the algebra homomorphism $U\g \to A$ to $\ol{U}=U\g'/(t^*)
\to A$ by $t\mapsto s+\rho_r$. Thus, $V$ can be considered as a
$\g'$--module on which $t^*$ acts trivially and $t$ by $s + \rho_r$.
Set
\[M_{\pm}=\ind^{\g'}_{(\g_\pm\rtimes\C t)\oplus\C t^*}\C_\pm\]
where $\g_\pm\rtimes\C t$ acts on the one--dimensional module
$\C_\pm$ by $0$ and $t^*$ as multiplication by $\pm 1$. Frobenius
reciprocity yields an isomorphism
\[\Xi:\Hom_{\g_+\oplus\g_-\oplus\C t^*}\lp M_+,M_-^*\ctensor
V^{\otimes n}\rp\to V^{\otimes n}\]
where $\ctensor$ is the completed tensor product.

Consider now the following system of partial differential equations
for a function $\Psi(z_0,\ldots, z_n)$ with values in $\Hom_{\g_+\oplus\g_-\oplus\C t^*}\lp
M_+,M_-^*\wh{\otimes} V^{\otimes n}\rp$
\begin{equation}\label{app-eq: rational-KZ}
\pde{\Psi}{z_k}=\frac{\hbar}{2\pi\iota}\lp \sum_{j\neq k}
\frac{\Omega_{kj}'}{z_k-z_j}\rp \Psi
\end{equation}
where $\Omega'=\Omega+t\otimes t^*+t^*\otimes t$ is the Casimir
tensor of $\g'$. One readily checks that, for any $1\leq i,j\leq n$,
\[\Xi\,\Omega'_{ij}\,\Xi^{-1}=\Omega_{ij}\aand
  \Xi\,\Omega'_{0i}\,\Xi^{-1}=s^{(i)}-\sum_{\substack{1\leq k\leq n\\k\neq i}}r_{ki}\]
Coupled with the change of variables $z_i=e^{u_i}$, $i=1,\ldots,n$,
this yields the following
\begin{prop}\label{app-prop: trig-rational}
Under the Frobenius reciprocity isomorphism 
\[\Xi : \Hom_{\g_+\oplus\g_-\oplus\C t^*}\lp M_+, M_-^*\wh{\otimes}V^{\otimes n}\rp\isom V^{\otimes n}\]
the restriction of \eqref{app-eq: rational-KZ} to $z_0=0$ coincides
with \eqref{app-eq: trig-KZ}.
\end{prop}

\subsection{} 

Denote the \EK quantization of a \fd Lie bialgebra $\ll$ by $\Uh\ll$.
By functoriality of quantization,
\[U_\hbar((\g_\pm\rtimes\C t)\oplus\C t^*)\cong
(\Uhg_\pm\rtimes\C[t])\otimes\C[t^*]\]
Set $M_{\pm}^q=\ind^{\Uhg'}_{U_\hbar((\g_\pm\rtimes\C t)\oplus\C t^*)}
\C_\pm$, where $U_\hbar(\g_\pm\rtimes\C t)$ acts trivially on $\C_\pm
\cong\C$ and $t^*$ as multiplication by $\pm 1$.

Regard $\V=V\fml$ as a $\Uhg'$--module via the homomorphism
$\Uhg'\cong U\g'\fml\to U\g\fml\to A\fml$, where the intermediate
map $\g'\to\g$ is given by $t^*\to 0$ and $t\to s+\rho_r$. Let $\wtR'$
be the $R$--matrix of $\Uhg'$. The following is a consequence of
the \KD theorem for $\g'$ \cite{etingof-kazhdan-quantization-1},
together with Proposition \ref{app-sec: proof-begin}. 

\begin{prop}\label{app-prop: ek-monodromy}
The monodromy representation \eqref{app-eq: monodromy}
is equivalent to the representation of $\Pi_n$ on $\Hom_{\Uh
(\g_+\oplus\g_-\oplus\C t^*)}\lp M_+^q,(M^q_-)^*\ctensor\V^
{\otimes n}\rp$ given by
\begin{align*}
b_i 	&\mapsto (i\, i+1) \wtR'_{i\,i+1} \\
\X_1 &\mapsto \wtR'_{10}\wtR'_{01}
\end{align*}
where $(M^q_-)^*$ is the right dual to $M^q_-$.\footnote
{the action of $a\in\Uhg'$ on $\phi\in(M^q_-)^*$ is given by
$a\phi=\phi\circ(S^{-1}a)$.}
\end{prop}

\subsection{}\label{app-sec: proof-end}

Since $\wtR'_{i\,i+1}$ acts on $\V^{\otimes n}$ as $\wtR_{i\,i+1}$,
Proposition \ref{app-prop: ek-monodromy} reduces the proof of
Theorem \ref{app-thm: etingof-geer} to computing the action of
$\Xi^q\,\wtR'_{10}\wtR_{01}'\,(\Xi^q)^{-1}$ on $V^{\otimes n}$,
where $\Xi^q$ is the isomorphism given by the composition
\begin{equation}\label{app-eq: big-diagram}
\xymatrix{
\Hom_{\Uh(\g_+\oplus\g_-\oplus\C t^*)}\lp M_+^q,(M^q_-)^*\ctensor \V^{\otimes n}
\rp \ar[dd]_{\Xi^q} \ar[rrr]^{\ev(\id_+^q)} &&& \lp(M^q_-)^*\ctensor 
\V^{\otimes n}\rp^{\Uhg_+} \ar[dd]^{\phi} \\
&&&\\
\V^{\otimes n} &&& \Hom_{\Uhg_+}\lp M_-^q, 
\V^{\otimes n}\rp \ar[lll]^{\ev(\id^q_-)}}
\end{equation}
and $\phi$ is the restriction of the natural identification $(M_-^q)^*
\otimes \V^{\otimes n} \cong \Hom\lp M_-^q, \V^{\otimes n}\rp$
to the subspace of $\Uhg_+$--invariant vectors.\footnote{
This is the reason for considering the right dual of $M^q_-$
instead of the left one as in \cite{etingof-geer}. For the latter,
the natural identification $\lp M_-^q\rp^*\otimes \V^{\otimes n} \cong
\Hom(M_-^q, \V^{\otimes n})$ does not restrict to the isomorphism
between the subspace of $\Uhg_+$--invariant and $\Uhg_+$--linear
morphisms.}

Write $\wtR'=\alpha_j\otimes \beta^j$, where $\{\alpha_j\}$ is
a basis of $\Uhg'_+$, $\{\beta^j\}$ is the dual basis of $\Uhg'_-$,
and the summation over $j$ is implicit. For a morphism
\[\Psi\in
\Hom_{\Uh(\g_+\oplus\g_-\oplus\C t^*)}\lp M_+^q,(M^q_-)^*\ctensor \V^{\otimes n}\rp\]
we compute $\X_1(\Xi^q(\Psi))$ in the following steps. In order to
make the computations more transparent we abusively assume
that $\Psi(\id^q_+)$ is an indecomposable tensor $m\otimes v_1
\otimes\cdots\otimes v_n$.

\begin{itemize}
\item[(a)] The computation below corrects equation (5.4) of \cite
{etingof-geer}.
\begin{align*}
\X_1(\Xi^q(\Psi)) &= \la \id^q_-, \wtR_{10}'\wtR_{01}'\Psi
\id^q_+\ra \\
&= \la \id^q_-, \beta^j\alpha_i\otimes \alpha_j\beta^i \otimes 
1^{\otimes n-1} (\Psi\id^q_+)\ra\\
&= \la \id_-^q, \beta^j\alpha_im\ra \alpha_j\beta^iv_1\otimes v_2\otimes
\cdots\otimes v_n\\
&= \la S^{-1}(\beta^j)\id_-^q, \alpha_im\ra \alpha_j\beta^iv_1\otimes v_2
\otimes\cdots\otimes v_n\\
&= \la \id_-^q, \alpha_im\ra e^{\hbar t}\beta^iv_1\otimes v_2\otimes\cdots
\otimes v_n \\
&= (e^{\hbar t})^{(1)}\,\Xi^q(\wtR_{01}\Psi)
\end{align*}
where we used the fact that $t^*$ acts trivially on $\V\ni v_1$ and the fifth
equality uses the fact that $\id^q_-$ is killed by $\Uhg_-$, that $t^*\id_-^q$
acts by $-1$ on $M_-^q$, and that the dual element to $(t^*)^k\in\Uhg'_-$
is $(\hbar t)^k/k!\in\Uhg'_+$.

\item[(b)] Write $\wtR=a_j\otimes b^j$, where $\{a_j\}$ is a basis of
$\Uhg_+$ and $\{b^j\}$ is the dual basis of $\Uhg_-$. Then
\begin{align*}
\X_1(\Xi^q(\Psi)) &= \la \id_-^q, a_im\ra e^{\hbar t}b^iv_1\otimes v_2
\otimes\cdots\otimes v_n \\
&= (e^{\hbar t}b^i)^{(1)}\la S^{-1}(a_i)\id_-^q, 
m\ra v_1\otimes v_2\otimes\cdots\otimes v_n \\
&= (e^{\hbar t}b^i)^{(1)} \la \id_-^q, m\ra
\Delta^{(n)}(S^{-1}(a_i)) \lp v_1\otimes v_2\otimes \cdots \otimes v_n\rp \\
&= (e^{\hbar t}b^i)^{(1)} \Delta^{(n)}(S^{-1}(a_i))
\Xi^q(\Psi)
\end{align*}
\item[(c)] Note that
\[\begin{split}
(b^i)^{(1)}\Delta^{(n)}(S^{-1}(a_i))
&=
m_{01}\left(\Id\otimes\Delta^{(n)}\circ\Id\otimes S^{-1}(\wtR_{21})\right)\\
&=
m_{01}\left(\Id\otimes\Delta^{(n)}\circ S\otimes\Id(\wtR_{21})\right)\\
&=
C\end{split}\]
where the second equality follows from the fact that $S\otimes S(\wtR)=\wtR$
and $C$ is the element defined in \eqref{app-eq: eg-element}). It follows that
$\X_1\in\Pi_n$ acts on $\V^{\otimes n}$ as
\[\X_1 \mapsto (e^{\hbar s}w)^{(1)} C\]
where $w=e^{\hbar \rho_r}$ under the identification 
of $\Uhg$ with $U\g[[\hbar]]$.
\item[(d)] To determine the element $w$, we restrict ourselves to the case
$s=0$ and $n=1$. In this case the monodromy is  trivial and hence we get
\[1=w\cdot m_{01}(S\otimes\Id(\wtR_{21}))=wu\]
where $u=S(b^i)a_i$ is the Drinfeld element. Thus $w=u^{-1}$, which
completes the proof of Theorem \ref{app-thm: etingof-geer}.
\end{itemize}

\section{Proof of Proposition \ref{prop: log-completion}}\label{app:J^n}

We shall prove the result for $\hloopsl{2}$. The corresponding
assertion for $\hloopgl{2}$ is proved similarly. The key step is
to draw the following consequence of \cite{guay-degeneration}

\begin{lem}\label{pf-lem: main}
The following elements are in $\J^n$ for any $n\geq 0$ and $l>0$
\begin{align*}
H_{0;n} 	&= H_0 + \frac{q-q^{-1}}{\hbar} \sum_{r=1}^n (-1)^r \cbin{n}{r} H_r \\
H_{l;n}	&= \frac{q-q^{-1}}{\hbar} \sum_{r=0}^n (-1)^r \cbin{n}{r}H_{l+r}
\end{align*}
\end{lem}

\noindent The proof of Lemma \ref{pf-lem: main} will be given in \S
\ref{ss:begin app proof}--\ref{ss:end pf lem}.

Let us prove that $\wt{H}_r\in \J^r$ using Lemma \ref{pf-lem: main}.
We will need the following easy

\begin{lem}
Let $\{X_k\}_{k\in\Z}$ be elements of a vector space $V$. For
each $t\in\Z$, and $0\leq m\leq n$, define
\[X_{n;t}^{(m)}=\sum_{s=0}^n (-1)^s \cbin{n}{s} s^m X_{s+t}\]
Then, for $m>0$
\[X_{n;t}^{(m)}=-n\sum_{r=0}^{m-1}\cbin{m-1}{r}X_{n-1;t+1}^{(r)}\]
In particular, $X_{n;t}^{(m)}$ can be written as a linear combination
of  $\{X_{n-k;t+k}^{(0)}\}_{1\leq k\leq m}$.
\end{lem}

In particular, taking $X_0=H_0$ and $X_k=\frac{q-q^{-1}}{\hbar}H_k$ for
$k\neq 0$, we see that $X_{n;t}^{(m)} \in \J^{n-m}$ since, by Lemma
\ref{pf-lem: main}, $X_{n-k;t+k}^{(0)}\in\J^{n-k}$. Since multiplication
by $\hbar$ maps $\J^n$ to $\J^{n+1}$, it follows that for any formal
power series $p(u) \in 1+ u\C[[u]]$ the following expression lies in $
\J^n$
\[X_{n;t}^{(p(u))}=\sum_{s=0}^n (-1)^s \cbin{n}{s} p(s\hbar)X_{s+t}\]

\noindent Taking $p(u)=\ds{\frac{u}{e^{u/2}-e^{-u/2}}}$, so that $\ds{
p(s\hbar)=\frac{\hbar}{q-q^{-1}}\frac{s}{[s]}}$, we see that
\begin{align*}
\wt{H}_r &= H_0 + \sum_{s=1}^r (-1)^s \cbin{r}{s}\frac{s}{[s]}H_s\\
&= H_0 + \frac{q-q^{-1}}{\hbar} \lp \sum_{s=1}^r (-1)^s \cbin{r}{s} p(s\hbar)
H_s\rp
\end{align*}
lies in $\J^r$, as claimed.

\subsection{Some notation}\label{ss:begin app proof}

For notational convenience, we set
$x_k = \frac{q-q^{-1}}{\hbar} H_k$ for $k\geq 1$ and $x_0=H_0$. 
Then we have
\[
\psi^+(z) = q^{H_0}\exp\lp \hbar \sum_{r\geq 1} x_rz^{-r}\rp
\]
and hence we obtain:
\begin{equation}\label{pf-eq: psi-explicit}
\psi_l = q^{H_0} \sum_{\lambda\vdash l} \hbar^{l(\lambda)} 
\frac{x_{\lambda}}{\prod_{i\geq 1}l_i!}
\end{equation}

Set $y_{l;n} = q^{-H_0}\frac{q-q^{-1}}{\hbar}\phi_{l;n}$. Then
using \eqref{pf-eq: psi-explicit} we have:

\begin{equation}\label{pf-eq: y-l}
y_{l;n} = \sum_{r=0}^n (-1)^r \cbin{n}{r} \lp \sum_{\lambda\vdash
l+r} \hbar^{l(\lambda)-1} \frac{x_{\lambda}}{\prod l_i!}\rp
\end{equation}

\begin{equation}\label{pf-eq: y-0}
y_{0;n} = \frac{1-e^{-\hbar x_0}}{\hbar} + \sum_{r=1}^n (-1)^r
\cbin{n}{r} \lp \sum_{\lambda\vdash r} \hbar^{l(\lambda)-1}
\frac{x_{\lambda}}{\prod l_i!}\rp
\end{equation}

\subsection{}

It is clear from the definitions that $y_{l;n}\in \J^n$ 
for every $l,n\geq 0$. We denote by $p_{l;n}^{(m)}$ the 
coefficient of $\frac{\hbar^{m-1}}{m!}$ in $y_{l;n}$. Thus we 
have the following expressions (here $l\geq 1$):
\begin{gather}
p_{l;n}^{(m)} = \sum_{r=0}^n (-1)^r \cbin{n}{r} \lp 
\sum_{\begin{subarray}{c}\lambda \vdash l+r \\ l(\lambda)=m 
\end{subarray}} \frac{m!x_{\lambda}}{\prod l_i!}\rp \label{pf-eq: p-l}\\
p_{0;n}^{(m)} = (-1)^{m-1}x_0^m + \sum_{r=1}^n
(-1)^r \cbin{n}{r} \lp \sum_{\begin{subarray}{c}\lambda\vdash r\\
l(\lambda)=m\end{subarray}} \frac{m!x_{\lambda}}{\prod l_i!}\rp
\label{pf-eq: p-0}
\end{gather}

We will prove the following stronger version of Lemma 
\ref{pf-lem: main}.

\begin{lem}\label{pf-lem: main2}
For each $l,n\geq 0$ and $m\geq 1$ we have:
\[
p_{l;n}^{(m)} \in \J^{n-m+1}
\]
\end{lem}

Note that the assertion of Lemma \ref{pf-lem: main} is $m=1$
case of that of Lemma \ref{pf-lem: main2}.

\subsection{Proof of Lemma \ref{pf-lem: main2}}

We begin by considering $l\geq 1$ case. In this case
we have the following relation for $m\geq 2$:

\begin{equation}\label{pf-eq: relation-l}
p_{l;n}^{(m)} = \sum_{t=1}^{l-1} p_{t;0}^{(m-1)}p_{l-t;n}^{(1)} - 
\sum_{k=0}^{n-1} p_{1;k}^{(1)}p_{l;n-k-1}^{(m-1)}
\end{equation}

We prove $p_{l;n}^{(m)} \in \J^{n-m+1}$ for every $l\geq 1$, 
$m\geq 1$, $n\geq 1$ by induction on $n$ and $l$ in the following
manner. Consider the base case of $n=1$:
\[
y_{l;1} = p_{l;1}^{(1)} + O(\hbar) \in \J
\]
which implies that $p_{l;1}^{(1)}\in \J$. For $n=1$ and $m\geq 2$
the statement is vacuous. Thus we have proved the assertion 
for $n=1$ and all $l,m\geq 1$.\\

Now we proceed to the induction step. Let us assume that 
$p_{l;n'}^{(m)} \in \J^{n'-m+1}$ for every $n'<n$ and $l,m\geq 1$.
Now the same assertion for $n'=n$ is proved for $m\geq 2$
by using \eqref{pf-eq: relation-l} and induction on $l$. The base 
case $l=1$ is established by \eqref{pf-eq: relation-l}:
\[
p_{1;n}^{(m)} = - \sum_{k=1}^{n-1} p_{1;k}^{(1)}p_{1;n-k-1}^{(m-1)}
\]
since all the terms on the \rhs have smaller $n$. Proceeding
by induction on $l$ we can prove the desired assertion for $n'=n$
 and for every $m\geq 2$, $l\geq 1$. The case $m=1$ follows 
 from the fact that
\[
y_{l;n} = p_{l;n}^{(1)} + \sum_{m\geq 2} \frac{\hbar^{m-1}}{m!}
p_{l;n}^{(m)} \in \J^n 
\]

\subsection{}

Next we consider the $l=0$ case. In this case
we will need the following relation (again for $m\geq 2$):

\begin{equation}\label{pf-eq: relation-0}
p_{0;n}^{(m)} = (-1)^{m-1}x_0^{m-1}p_{0;n-1}^{(1)}
-\sum_{k=0}^{n-2} p_{1;k}^{(1)}p_{0;n-1-k}^{(m-1)}
\end{equation}

Again $p_{0;n}^{(m)}\in \J^{n-m+1}$ is proved by an induction 
argument, similar to the one given above, using 
\eqref{pf-eq: relation-0}, combined with the result of the 
previous section.

\subsection{Proof of \eqref{pf-eq: relation-l} and \eqref{pf-eq: relation-0}}\label{ss:end pf lem}

The proofs of relations \eqref{pf-eq: relation-l} and 
\eqref{pf-eq: relation-0} are similar. We provide the main steps
in the proof of \eqref{pf-eq: relation-l} and leave a few 
straightforward checks to the reader.\\

For the proof, it will be convenient to write $p_{l;n}^{(m)}$ as:
\[
p_{l;n}^{(m)} = \sum_{r=0}^n (-1)^r \cbin{n}{r}\lp
\sum_{\begin{subarray}{c} a_1,\ldots,a_m\geq 1 \\
a_1+\cdots+a_m = r+l\end{subarray}}x_{a_1}\ldots x_{a_m}\rp
\]

which implies the following verification:

\begin{equation}\label{pf-eq: lhs}
\sum_{t=1}^{l-1} p_{t;0}^{(m-1)}p_{l-t;n}^{(1)} - p_{l;n}^{(m)}
= \sum_{s=0}^{n-1} (-1)^s \cbin{n}{s+1}\lp \sum_{
\begin{subarray}{c} a_1,\ldots,a_{m-1}\geq 1 \\ l\leq 
a_1+\cdots+a_{m-1}\leq l+s\end{subarray}}
x_{a_1}\cdots x_{a_{m-1}}\rp
\end{equation}

Let us denote the expression obtained above by $g_{m,l}(n)$.
Recall that the equation \eqref{pf-eq: relation-l} is equivalent to
\[
g_{m,l}(n) = \sum_{k=0}^{n-1}p_{1;k}^{(1)}p_{l;n-k-a}^{(m-1)}
\]

This equation can be verified in the following manner. It is 
easy to check that the claimed equation holds for $n=1$. 
Moreover both sides satisfy the following recurrence relation:
\[
F_{m,l}(n+1) - F_{m,l}(n) + F_{m,l+1}(n) = p_{1;n}^{(1)}p_{l;0}^{(m-1)}
\]
which implies the desired assertion by induction on $n$.


\providecommand{\bysame}{\leavevmode\hbox to3em{\hrulefill}\thinspace}
\providecommand{\MR}{\relax\ifhmode\unskip\space\fi MR }
\providecommand{\MRhref}[2]{%
  \href{http://www.ams.org/mathscinet-getitem?mr=#1}{#2}
}
\providecommand{\href}[2]{#2}

\end{document}